%From leonhard@math.rutgers.edu Fri Sep 15 16:57:08 1995
%Date: Fri, 15 Sep 95 10:54:48 EDT
%From: Alice Leonhardt <leonhard@math.rutgers.edu>
%To: shlhetal@MATH.HUJI.AC.IL
%Subject: sh481
%
%
%Andrzej:
%
%last one - 481 - alice
%------
\input amstex
\NoBlackBoxes
%\def\cite #1{\rm[#1]}
%%%\magnification=\magstep 1
\documentstyle {amsppt}
\topmatter
\title {Was Sierpinski Right III? \\
Can Continuum-c.c. times c.c.c. be continuum-c.c.? \\
Sh481} \endtitle
\rightheadtext{Was Sierpinski Right III}
\author {Saharon Shelah \thanks {\null\newline
Typed 5/92 - Latest Revision 1/13/95 \null\newline
I thank Alice Leonhardt for the beautiful typing \null\newline
\S1 corrected 4/94 \null\newline
Partially supported by the basic research fund, Israeli
Academy} \endthanks} \endauthor
%Previous version - 10/7/94
\affil {Institute of Mathematics \\
The Hebrew University \\
Jerusalem, Israel
\medskip
Rutgers University \\
Department of Mathematics \\
New Brunswick, NJ  USA} \endaffil

\abstract  We prove the consistency of: if  $B_1$, $B_2$ are Boolean algebras 
satisfying the c.c.c. and the $2^{\aleph_0}$-c.c. respectively then  $B_1 
\times  B_2$ satisfies the $2^{\aleph_0}$-c.c. 

We start with a universe with a Ramsey cardinal (less suffice).
\endabstract 
\endtopmatter
\document
\expandafter\ifx\csname bib4plain.tex\endcsname\relax
  \expandafter\gdef\csname bib4plain.tex\endcsname{}
\else \message{Hey!  Apparently you were trying to \string twice.   This does not make sense.}
\errmessage{Please edit your file (probably \jobname.tex) and remove
any duplicate ``\string\input'' lines} \fi

%  This file should be inputted if you want to use 
%  bibtex fom within plain TeX. 
      % Not really need for standard
       % bibtex files, but these commands
\def\renewcommand{\newcommand}	       % are used in our literal-unsrt.bst
\edef\cite{\the\catcode`@}%
\catcode`@ = 11
\let\@oldatcatcode = \cite
\chardef\@letter = 11
\chardef\@other = 12
%
%
% Next come some things that will be useful later.
%
% Make an outer definition into an inner one (due to Chris Thompson).
% The arguments should be the control sequence to be defined, and the
% new of the \outer control sequence, as characters; the control
% sequence #1 is defined to be just the same as \csname#2\endcsname, but
% not \outer.  For example, \@innerdef\innernewcount{newcount} would
% define \innernewcount to be a non-outer version of \newcount.
%
\def\@innerdef#1#2{\edef#1{\expandafter\noexpand\csname #2\endcsname}}%
%
% We use \@innerdef to make some of our allocations local, because
% Eplain includes our code inside a conditional.  We put @'s in the
% names to minimize the (already small) chance of conflicts.
%
\@innerdef\@innernewcount{newcount}%
\@innerdef\@innernewdimen{newdimen}%
\@innerdef\@innernewif{newif}%
\@innerdef\@innernewwrite{newwrite}%
%
%
% Swallow one parameter.
%
\def\@gobble#1{}%
%
%
% Use TeX 3.0's \inputlineno to get the line number, for better error
% messages, but if we're using an old version of TeX, don't do anything.
%
\ifx\inputlineno\@undefined
   \let\@linenumber = \empty % Pre-3.0.
\else
   \def\@linenumber{\the\inputlineno:\space}%
\fi
%
%
% The following macro \@futurenonspacelet (from the TeXbook) behaves
% essentially like \futurelet except that it discards any implicit or
% explicit space tokens that intervene before a nonspace is scanned:
%
\def\@futurenonspacelet#1{\def\cs{#1}%
   \afterassignment\@stepone\let\@nexttoken=
}%
\begingroup % The grouping here avoids stepping on an outside use of `\\'.
\def\\{\global\let\@stoken= }%
\\ % now \@stoken is a space token (\\ is a control symbol, so that
   % space after it is seen).
\endgroup
\def\@stepone{\expandafter\futurelet\cs\@steptwo}%
\def\@steptwo{\expandafter\ifx\cs\@stoken\let\@@next=\@stepthree
   \else\let\@@next=\@nexttoken\fi \@@next}%
\def\@stepthree{\afterassignment\@stepone\let\@@next= }%
%
%
% \@getoptionalarg\CS gets an optional argument from the input, enclosed
% in brackets, then expands \CS.  We set \@optionalarg to \empty if we
% don't find one, otherwise to the text of the argument.  This assumes
% the brackets don't have a funny category code.
%
\def\@getoptionalarg#1{%
   \let\@optionaltemp = #1%
   \let\@optionalnext = \relax
   \@futurenonspacelet\@optionalnext\@bracketcheck
}%
%
% The \expandafter's in this macro let us avoid the use of \aftergroup,
% which is somewhat more expensive.
%
\def\@bracketcheck{%
   \ifx [\@optionalnext
      \expandafter\@@getoptionalarg
   \else
      \let\@optionalarg = \empty
      % We can't do the \temp after the \fi, because then the \temp gets
      % in the way of reading the optional argument from the input, if
      % we do have one.
      \expandafter\@optionaltemp
   \fi
}%
\def\@@getoptionalarg[#1]{%
   \def\@optionalarg{#1}%
   \@optionaltemp
}%
%
%
% From LaTeX.
%
\def\@nnil{\@nil}%
\def\@fornoop#1\@@#2#3{}%
\def\@for#1:=#2\do#3{%
   \edef\@fortmp{#2}%
   \ifx\@fortmp\empty \else
      \expandafter\@forloop#2,\@nil,\@nil\@@#1{#3}%
   \fi
}%
\def\@forloop#1,#2,#3\@@#4#5{\def#4{#1}\ifx #4\@nnil \else
       #5\def#4{#2}\ifx #4\@nnil \else#5\@iforloop #3\@@#4{#5}\fi\fi
}%
\def\@iforloop#1,#2\@@#3#4{\def#3{#1}\ifx #3\@nnil
       \let\@nextwhile=\@fornoop \else
      #4\relax\let\@nextwhile=\@iforloop\fi\@nextwhile#2\@@#3{#4}%
}%
%
%
% This macro tests if a file \jobname.#1 exists, and sets \if@fileexists
% appropriately.  If an optional argument is given, it is used as the
% root part of the filename instead of \jobname.
%
\@innernewif\if@fileexists
\def\@testfileexistence{\@getoptionalarg\@finishtestfileexistence}%
\def\@finishtestfileexistence#1{%
   \begingroup
      \def\extension{#1}%
      \immediate\openin0 =
         \ifx\@optionalarg\empty\jobname\else\@optionalarg\fi
         \ifx\extension\empty \else .#1\fi
         \space
      \ifeof 0
         \global\@fileexistsfalse
      \else
         \global\@fileexiststrue
      \fi
      \immediate\closein0
   \endgroup
}%
%
%
%% [[[start of BibTeX-specific stuff]]]
%
% Now come the four main LaTeX commands and their associated .aux
% commands.  Just as in LaTeX, \bibliographystyle defines the BibTeX
% style name (.bst file, that is), and \bibliography defines the
% database (.bib) file(s).  The corresponding .aux-file commands are
% \bibstyle and \bibdata, which are there only for BibTeX's (but not
% LaTeX's) use.
%
\def\bibliographystyle#1{%
   \@readauxfile
   \@writeaux{\string\bibstyle{#1}}%
}%
\let\bibstyle = \@gobble
%
% As well as writing the \bibdata command to tell BibTeX which .bib
% files to read, we read the .bbl file that BibTeX (or a person,
% conceivably) has produced.  We use \bblfilebasename as the root of the
% filename to read; this defaults to \jobname.
%
\let\bblfilebasename = \jobname
\def\bibliography#1{%
   \@readauxfile
   \@writeaux{\string\bibdata{#1}}%
   \@testfileexistence[\bblfilebasename]{bbl}%
   \if@fileexists
      % We just output a non-discardable item (the `whatsit' with the
      % \bibdata command).  This means that the glue that will be
      % inserted next (\parskip or \baselineskip, most likely) will be a
      % legal breakpoint.  Most likely, this is after some kind of
      % heading, where we don't want to allow a page break.  So:
      \nobreak
      \@readbblfile
   \fi
}%
\let\bibdata = \@gobble
%
% The \nocite{label,label,...} command writes its argument to \@auxfile,
% unless instructed not to, but produces no text in the document.  Both
% the \nocite and \cite commands produce \citation commands in the .aux file.
%
\def\nocite#1{%
   \@readauxfile
   \@writeaux{\string\citation{#1}}%
}%
\@innernewif\if@notfirstcitation
%
% \cite[note]{label,label,...} produces the citations for the labels as
% well.  If the optional argument `note' is present, it's added after
% the labels.  Since \cite calls \nocite to do its .aux-file writing,
% \cite doesn't need to call \@readauxfile (\nocite does).
%
\def\cite{\@getoptionalarg\@cite}%
%
% Typeset the citations for the labels in #1, followed by the note, if
% it exists.  To change the citation's format in the text, redefine one
% or more `\print...' macros, whose defaults appear later in this file.
%
\def\@cite#1{%
   % Remember the optional argument, in case one of the macros we call
   % below ends up looking for an optional argument itself.  For
   % example, if a \cite[note] triggers reading the .aux file, then the
   % [note] would be clobbered, since \@testfileexistence looks for an
   % optional arg.
   \let\@citenotetext = \@optionalarg
   % Start printing the text, beginning with a left bracket by default.
   \printcitestart
   % It's complicated, but because \nocite puts a `whatsit' onto the list,
   % \nocite should follow \printcitestart.  It's conceivable, but very
   % unlikely, that this `whatsit' will cause a problem (glue that doesn't
   % disappear when you want it to is the most likely symptom), requiring
   % a change either to \printcitestart or to the label that the .bst file
   % produces.
   \nocite{#1}%
   \@notfirstcitationfalse
   \@for \@citation :=#1\do
   {%
      \expandafter\@onecitation\@citation\@@
   }%
   \ifx\empty\@citenotetext\else
      \printcitenote{\@citenotetext}%
   \fi
   \printcitefinish
}%
\def\@onecitation#1\@@{%
   \if@notfirstcitation
      \printbetweencitations
   \fi
   \expandafter \ifx \csname\@citelabel{#1}\endcsname \relax
      \if@citewarning
         \message{\@linenumber Undefined citation `#1'.}%
      \fi
      % Give it a dummy definition:
      \expandafter\gdef\csname\@citelabel{#1}\endcsname{%
% Change: marginal remark added, goldstrn@math.huji.ac.il, 
% goldstern@tuwien.ac.at, May 1996 mg
%  !!! change !!!
\strut
\vadjust{\vskip-\dp\strutbox
\vbox to 0pt{\vss\parindent0cm \leftskip=\hsize 
\advance\leftskip3mm
\advance\hsize 4cm\strut\openup-4pt 
\rightskip 0cm plus 1cm minus 0.5cm ?  #1 ?\strut}}
         {\tt
            \escapechar = -1
            \nobreak\hskip0pt
            \expandafter\string\csname#1\endcsname
            \nobreak\hskip0pt
         }%
      }%
   \fi
   % Now produce the text, whether it was undefined or not.
   \csname\@citelabel{#1}\endcsname
   \@notfirstcitationtrue
}%
%
% Given a label `foo', the macro `\b@foo' is supposed to
% hold the text that should be produced.
%
\def\@citelabel#1{b@#1}%
%
% So, how does a citation label get defined?  When we read the .bbl file
% (below), a \bibitem writes out a \@citedef command.  And when we read
% the \@citedef, we define \@citelabel{#1}, where #1 is the user's
% label.
%
\def\@citedef#1#2{\expandafter\gdef\csname\@citelabel{#1}\endcsname{#2}}%
%
%
% Reading the .bbl file also produces the typeset bibliography.  Please
% notice, however, that we do not produce the title for the references
% (e.g., `References'), as LaTeX does.  The formatting and spacing of
% that title, whether it should go into the headline, and so on, are all
% things determined by your format.  We cannot know those things in
% advance.  If you wish, you can define \bblhook to produce the title.
% Or just do it before the \bibliography command.
%
\def\@readbblfile{%
   % Define a counter to tell us which item number we are on, unless
   % we've already defined it (because the document has more than one
   % bibliography).
   \ifx\@itemnum\@undefined
      \@innernewcount\@itemnum
   \fi
   \begingroup
      \def\begin##1##2{%
         % ##1 is just `thebibliography'.
         % ##2 is the widest label.
         % We set (new dimen) \biblabelwidth based on the widest label
         \setbox0 = \hbox{\biblabelcontents{##2}}%
         \biblabelwidth = \wd0
      }%
      \def\end##1{}% ##1 is `thebibliography' again.
      %
      % Here we have two possibilities:
      % \bibitem[typesetlabel]{citationlabel}
      % \bibitem{citationlabel}
      % If we have the second of these, the citations are numbered, starting
      % from one; we use our own count register \@itemnum for this.
      %
      \@itemnum = 0
      \def\bibitem{\@getoptionalarg\@bibitem}%
      \def\@bibitem{%
         \ifx\@optionalarg\empty
            \expandafter\@numberedbibitem
         \else
            \expandafter\@alphabibitem
         \fi
      }%
      \def\@alphabibitem##1{%
         % Need \xdef here for various reasons.
         \expandafter \xdef\csname\@citelabel{##1}\endcsname {\@optionalarg}%
         % Left-justify alpha labels, unless \biblabel{pre,post}contents
         % are already defined.
         \ifx\biblabelprecontents\@undefined
            \let\biblabelprecontents = \relax
         \fi
         \ifx\biblabelpostcontents\@undefined
            \let\biblabelpostcontents = \hss
         \fi
         \@finishbibitem{##1}%
      }%
      \def\@numberedbibitem##1{%
         \advance\@itemnum by 1
         \expandafter \xdef\csname\@citelabel{##1}\endcsname{\number\@itemnum}%
         % Right-justify numeric labels, unless \biblabel{pre,post}contents
         % are already defined.
         \ifx\biblabelprecontents\@undefined
            \let\biblabelprecontents = \hss
         \fi
         \ifx\biblabelpostcontents\@undefined
            \let\biblabelpostcontents = \relax
         \fi
         \@finishbibitem{##1}%
      }%
      \def\@finishbibitem##1{%
         \biblabelprint{\csname\@citelabel{##1}\endcsname}%
         \@writeaux{\string\@citedef{##1}{\csname\@citelabel{##1}\endcsname}}%
         \ignorespaces
      }%
      %
      % Do the printing (we're producing the bibliography, remember).
      %
      \let\em = \bblem
      \let\newblock = \bblnewblock
      \let\sc = \bblsc
      % Punctuation won't affect spacing;
      \frenchspacing
      % the penalties below are from LaTeX's [article,book,report].sty;
      \clubpenalty = 4000 \widowpenalty = 4000
      % the next two values come from LaTeX's \sloppy command;
      \tolerance = 10000 \hfuzz = .5pt
      \everypar = {\hangindent = \biblabelwidth
                      \advance\hangindent by \biblabelextraspace}%
      \bblrm
      % the \parskip is a guess at what looks good;
      \parskip = 1.5ex plus .5ex minus .5ex
      % and the space between label and text comes from LaTeX's \labelsep.
      \biblabelextraspace = .5em
      \bblhook
      \input \bblfilebasename.bbl
   \endgroup
}%
%
% The widest label's width is useful for redefining \biblabelprint;
% you redefine \biblabelwidth, in effect, by redefining the
% \biblabelcontents macro that appears below.  And \biblabelextraspace,
% which is redefinable inside \bblhook, is added to \biblabelwidth to
% determine the amount of hanging indentation.
%
\@innernewdimen\biblabelwidth
\@innernewdimen\biblabelextraspace
%
% Now come the main macros that are related to the printing of the
% bibliography.  Since you might want to redefine them, they are given
% default definitions outside of \@readbblfile.
%
% The first one controls the printing of a bibliography entry's label.
% If you change it, make sure that it starts with something like
% \noindent or \indent or \leavevmode that puts TeX into horizontal mode
% (even if the label itself is empty); otherwise, the hanging
% indentation will get messed up in certain circumstances.
%
\def\biblabelprint#1{%
   \noindent
   \hbox to \biblabelwidth{%
      \biblabelprecontents
      \biblabelcontents{#1}%
      \biblabelpostcontents
   }%
   \kern\biblabelextraspace
}%
%
% If you are using numeric labels, and you want them left-justified
% (numeric labels by default are right-justified), do something like:
%     \def\biblabelprecontents{\relax}
%     \def\biblabelpostcontents{\hss}
%
% By default the labels are typeset in \bblrm, and enclosed in brackets.
%
\def\biblabelcontents#1{{\bblrm [#1]}}%
%
% The main text, too, is typeset using \bblrm, which is \rm by default.
%
\def\bblrm{\rm}%
%
% Emphasis for producing, e.g., titles, is done with \it by default.
%
\def\bblem{\it}%
%
% Some styles use a caps-and-small-caps font for author names.  LaTeX
% defines an \sc command but plain TeX doesn't, so we need one here.
% The definition below doesn't load the font unless it's needed, but it
% tries to load only the 10pt version, because it might not exist at
% other point sizes.
%
\def\bblsc{\ifx\@scfont\@undefined
              \font\@scfont = cmcsc10
           \fi
           \@scfont
}%
%
% The major parts of an entry are separated with \bblnewblock.  The
% numbers below are taken from LaTeX's `article' style.
%
\def\bblnewblock{\hskip .11em plus .33em minus .07em }%
%
% Here's where you stick any other bibliography-formatting goodies, or
% redefine the values above.
%
\let\bblhook = \empty
%
%
% Here are the four default definitions for formatting the in-text
% citations.  These are what you redefine (after your \input btxmac but
% before your \bibliography) to get parens instead of brackets, or
% superscripts, or footnotes, or whatever.
%
\def\printcitestart{[}%         left bracket
\def\printcitefinish{]}%        right bracket
\def\printbetweencitations{, }% comma, space
\def\printcitenote#1{, #1}%     comma, space, note (if it exists)
%
% That scheme is pretty flexible.  For example you could use
%     \def\printcitestart{\unskip $^\bgroup}
%     \def\printcitefinish{\egroup$}
%     \def\printbetweencitations{,}
%     \def\printcitenote#1{\hbox{\sevenrm\space (#1)}}
%     \font\eighttt = cmtt8
%     \scriptfont\ttfam = \eighttt
% to get superscripted in-text citations.  (The scriptfont stuff
% exists only to print an undefined citation; it's in cmtt8 because
% there is no cmtt7.)  To get something radically different, however,
% you'll have to define your own \cite command.
%
% When we read `\citation' from the .aux file, it means nothing.
%
\let\citation = \@gobble
%
%
% Now comes the stuff for dealing with LaTeX's \newcommand.  As
% mentioned earlier, this \newcommand will redefine a preexisting
% command; that's different from how LaTeX's \newcommand behaves.
%
\@innernewcount\@numparams
%
% \newcommand{\foo}[n]{text} defines the control sequence \foo to have
% n parameters, and replacement text `text'.
%
\def\newcommand#1{%
   \def\@commandname{#1}%
   \@getoptionalarg\@continuenewcommand
}%
%
% Figure out if this definition has parameters.
%
\def\@continuenewcommand{%
   % If no optional argument, we have zero parameters.  Otherwise, we
   % have that many.
   \@numparams = \ifx\@optionalarg\empty 0\else\@optionalarg \fi \relax
   \@newcommand
}%
%
% \@numparams is how many arguments this command has.  The name of the
% command is \@commandname.  The replacement text for the new macro is #1.
%
\def\@newcommand#1{%
   \def\@startdef{\expandafter\edef\@commandname}%
   \ifnum\@numparams=0
      \let\@paramdef = \empty
   \else
      \ifnum\@numparams>9
         \errmessage{\the\@numparams\space is too many parameters}%
      \else
         \ifnum\@numparams<0
            \errmessage{\the\@numparams\space is too few parameters}%
         \else
            \edef\@paramdef{%
               % This is disgusting, but \loop doesn't work inside \edef,
               % because \body isn't defined.
               \ifcase\@numparams
                  \empty  No arguments.
               \or ####1%
               \or ####1####2%
               \or ####1####2####3%
               \or ####1####2####3####4%
               \or ####1####2####3####4####5%
               \or ####1####2####3####4####5####6%
               \or ####1####2####3####4####5####6####7%
               \or ####1####2####3####4####5####6####7####8%
               \or ####1####2####3####4####5####6####7####8####9%
               \fi
            }%
         \fi
      \fi
   \fi
   \expandafter\@startdef\@paramdef{#1}%
}%
%
%% [[[end of BibTeX-specific stuff]]]
%
%
% Names of references (arguments given in the \cite and \nocite
% commands) and file names (arguments given in the \bibliography and
% \bibliographystyle commands) are recorded in \jobname.aux, called the
% \@auxfile in these macros.  Here's how they get read in.
%
\def\@readauxfile{%
   \if@auxfiledone \else % remember: \@auxfiledonetrue if \noauxfile is defined
      \global\@auxfiledonetrue
      \@testfileexistence{aux}%
      \if@fileexists
         \begingroup
            % Because we might be in horizontal mode when \@readauxfile
            % is called (if it's in response to a \cite or \nocite), we
            % want to ignore all the would-be spaces at the ends of
            % lines in the aux file.  Fortunately, it's highly unlikely
            % an end-of-line might actually be desired.
            % And because we don't change the category code of anything
            % but @, primitives like \gdef can't be used to define labels
            % in the aux file.  The solution adopted by btxmac.tex is to
            % write `\@citedef{LABEL}{DEFINITION}' to the aux file, and
            % use \csname on LABEL.
            \endlinechar = -1
            \catcode`@ = 11
            \input \jobname.aux
         \endgroup
      \else
         \message{\@undefinedmessage}%
         \global\@citewarningfalse
      \fi
      \immediate\openout\@auxfile = \jobname.aux
   \fi
}%
%
% The \@readauxfile macro does all that work the first time it's called.
% Since it's called once for every \cite, \nocite, \bibliography, and
% \bibliographystyle command that the user issues, we need to remember
% whether the work's been done.  It's considered done if we're not to do
% it---that is, if \noauxfile is defined.
%
\newif\if@auxfiledone
\ifx\noauxfile\@undefined \else \@auxfiledonetrue\fi
%
% It's conceivable you'd want to change how other characters are read;
% to do that, change their category code before doing \input btxmac.
%
%
% After reading the .aux file, \@readauxfile opens it for writing.
% The \@writeaux macro does the actual writing (as long as
% \noauxfile is undefined).
%
\@innernewwrite\@auxfile
\def\@writeaux#1{\ifx\noauxfile\@undefined \write\@auxfile{#1}\fi}%
%
%
% A macro package that uses btxmac.tex might define
% \@undefinedmessage (before doing an \input btxmac).
%
\ifx\@undefinedmessage\@undefined
   \def\@undefinedmessage{No .aux file; I won't give you warnings about
                          undefined citations.}%
\fi
%
% Even if citations are undefined, we want to complain only if
% \@citewarningtrue.  The default is to set \@citewarningtrue unless
% \noauxfile is defined.  Again, a macro package that uses
% btxmac.tex might want to redefine this.
%
\@innernewif\if@citewarning
\ifx\noauxfile\@undefined \@citewarningtrue\fi
%
%
% Finally, before leaving we restore @'s old category code.
%
\catcode`@ = \@oldatcatcode

\newpage

\head {\S0 Introduction} \endhead
\bigskip

We heard the problem from Velickovic who got it from Todorcevic,
it says ``are there  $P$,  a c.c.c. forcing notion, and  $Q$  is a
$2^{\aleph_0}$-c.c. forcing such that  $P \times Q$ is not
$2^{\aleph_0}$-c.c.?"  We can phrase 
it as a problem of cellularity of Boolean algebras or topological spaces.
\smallskip
\noindent
We give a negative answer even for  $2^{\aleph_0}$ regular, this
by proving the consistency of the negation.  The proof is close to
\cite{Sh:288},\S3 which continues \cite{Sh:276},\S2 and is close to 
\cite{Sh:289}.  A recent use is \cite{Sh:473}.

We start with $V \models ``\lambda$ is a Ramsey cardinal", use c.c.c. forcing
blowing the continuum to $\lambda$.  Originally the paper contained the
consistency of e.g. $2^{\aleph_0} \rightarrow [\aleph_2]^2_3,2^{\aleph_0}$
the first $k^2_2$-Mahlo, (weakly inaccessible)(remember $k^2_2 < \omega$) but
the theorem presented arrive here to satisfactory state (for me) 
earlier.  See more \cite{Sh:546}.  I thank Mariusz Rabus for corrections.  
\bigskip

\centerline {$* \qquad * \qquad *$}
\bigskip

\noindent
What problems do \cite{Sh:276}, \cite{Sh:288}, \cite{Sh:289}, \cite{Sh:473} and
\cite{Sh:481} raise?  The most important are (we state the simplest uncovered
case for each point):
\bigskip

\subhead {A Question} \endsubhead  1) Can we get e.g. $CON(2^{\aleph_0} 
\rightarrow [\aleph_2]^2_3)$; more generally raise $\mu^+$ to higher.
\newline
2)  Can we get $CON(\aleph_\omega > 2^{\aleph_0} \rightarrow [\aleph_1]^2_3)$;
generally lower $2^\mu$, the exact $\aleph_n$ seems to me less
exciting. \newline
3)  Can get e.g. $CON(2^\mu > \lambda \rightarrow [\mu^+]^2_3)$? 
\bigskip

\noindent
Also concerning \cite{Sh:473}.
\subhead{B Question} \endsubhead  1) Can we get the continuity on a non-
meagre set for functions $f:{}^\kappa 2 \rightarrow {}^\kappa 2$? \newline
2)  what can we say on continuity of $2$-place functions? \newline
3)  What about $n$-place functions (after \cite{Sh:288}).
\bigskip

\subhead {C Question} \endsubhead  1)  Can we get e.g. for $\mu = \mu^{< \mu}
> \aleph_0$,$CON$(if $P$ is $2^\mu$-c.c.,$Q$ is $\mu^+$-c.c. then
$P \times Q$ is $2^\mu$-c.c)? \newline
2)  Can we get e.g. $CON$ (if $P$ is $2^{\aleph_0}$-c.c.,$Q$ is 
$\aleph_2$-c.c. then $P \times Q$ is $2^{\aleph_0}$-c.c.) \newline
3)  Can we get e.g. $CON(2^{\aleph_0} > \lambda > \aleph_0$, and if 
$P$ is $\lambda$-c.c.,$Q$ is $\aleph_2$-c.c. then $P \times Q$ is
$\lambda$-c.c.) \newline
On A1 see \cite{Sh:546}.
\medskip

\noindent
\underbar{Discussion}  Maybe the solution to (A1) is by using squared 
demand and if in $\delta < \lambda,
cf(\delta) > \mu$ we guess $\langle N_s:s \in [B]^{\le 2} \rangle,
{\underset\tilde {}\to c}$, try to by $Q_\delta$ to add a large subset of
$B$ on which only two colours appear; but we want to do it also when
$\text{otp}(B) > \mu^+$.  Naturally we assume that if $\delta' < \delta,
cf(\delta') > \mu^+,\delta' = \text{ Sup}(B \cap \delta')$ this was done
$\langle N_s:s \in [B \cap \delta']^{\le 2} \rangle$, but we need more:
including dividing $B$ to $\mu$ set on each only two colours (by $Q_\delta$).
To do this and have $\lambda,k^2_3$-Mahlo (rather than measurable or
$(\lambda \rightarrow^+ (\omega_1)^{<\omega}_2)$ we have to use a very strong
diamond.
\medskip

For problem (A2) the natural thing is to use systems $\bar N = \langle N_s:
s \in [B]^{\le 2} \rangle$ which are not end extension systems.  Then it is
natural to use the forcing on this stronger; ``specializing" not only the
colouring but all $P_{N_s \cap \lambda}$-names of ordinals (as defined in
\S8).  This required a suitable squared diamond; this has not yet been
clarified (actually we need somewhat less than $\bar N,S$.
\medskip

But for problem (A3) a weaker version of this suggest itself.  As in the
solution of 1, $\lambda$ is $k^2_2$-Mahlo $\left< \langle N^\delta_s:s \in
[B^\delta]^{\le 2} \rangle:\delta \in S \right>$ is such that $N^\delta_s
\subseteq H(\delta^+)$, (we think of $N^\delta_s$ as guessing the isomorphism
type over $H(\delta)$.  Now we have to define a preliminary forcing $R$,
$\lambda$-complete or at least strategically $\lambda$-complete, satisfying
the $\lambda^+$-c.c.  So we have ``copies" of $\langle N^\delta_s:s \in
[B^\delta]^{\aleph 2} \rangle$ which behave like $\triangle$-systems.
[Saharon].
\medskip

But if want to get tree like systems (ease requirement on forcing) we need
more (enough dependency).  For simplicity $\lambda = cf(\chi) = \chi^\lambda$
and use the following instead forcing and do it with. \newline
We can have in $V$ (or force), \newline
$\bar C = \langle C_\delta:\delta \in S \rangle$ a square, $S \subseteq \chi,
\bar S =: \{ \delta \in S:cf(\delta) = \lambda \} \subseteq \chi$
stationary, $[\delta \in S \Rightarrow \text{ otp}(C_\delta) \le \lambda)$,
we have squared diamond $\bar C = \langle C_\delta:\delta \in S \rangle$, and
we choose for $\delta \in S$,$B_\delta \prec {\frak A}_\delta,\| B_\delta \|
\le |\omega + \text{ otp}(C_\delta)|$,$\delta_1 \in C_{\delta_2}
\Rightarrow B_{\delta_1} \prec B_{\delta_2} \and \langle B_\delta:\delta \in
C_{\delta} \cup \{ \delta_1\} \rangle \in B_{\delta_1}$,
$\delta_2 = \text{ sup }C_{\delta_2} = B_{\delta_1} = \dsize \bigcup
_{\delta \in C_{\delta_1}} B^*_\delta$. \newline
Now we can copy the squared diamond
$\left< \langle N_{\delta,s}:s \in [B_\delta]^{\le 2} \rangle:\delta <
\lambda \right>$ \newline
getting $\left< \langle N^*_{\delta,s}:s \in [B^*_\delta]^{\le 2} \rangle:
\delta \in S \right>$.  We then define $\langle P_i,Q_i,A_i:i < \chi
\rangle,|a_i| \le \chi_2$ (or $|a_i| < \kappa$).
\medskip

Concerning (B1) the expected theorem holds.  For $2$-place function, note
that the Sierpinski colouring can be viewed as a function from ${}^\mu 2$
to $\{ 0_\mu,1_\mu \} \subseteq {}^\mu 2$.  So the \underbar{best} we can
hope for is

(B2)$'$  can we get the consistency of $(*)_\mu$ for any $2$-place function
$f$ from ${}^\mu 2$ to ${}^\mu 2$ there are (everywhere) non-meagre
$A \subseteq {}^\mu 2$ and continuous functions $f_0,f_2$ from $A$ to
${}^\mu 2$ such that $(\forall \eta,\nu \in A)[f(\eta,\nu) \in
\{ f_0(\eta,\nu),f_1(\eta,nu)\}]$.
\medskip

So we have to put together the proofs of \cite{Sh:473} 
(continuity on non-meagre),
\cite{Sh:288} ($2^{\aleph_0} \rightarrow [\aleph_1]^2_3$), using the $k^2_2$-
Mahlo only and replace $\aleph_0$ by $\mu$.
\medskip

For problem (B3), $\mu = \aleph_0$ we have to generalize \cite{Sh:288},\S3.  
But
also for $\mu > \aleph_0$, we have to consider what can be said on the
partition of trees (see \cite{Sh:288},\S4 for a positive answer for 
large cardinal (indestructible measurable $n^* < \omega$).
\medskip

Concerning (C2), the problem with the approach to (A1) is ``why should
$Q_\delta$ from \cite{Sh:481},1.7 satisfies $Q^2_\delta$ is c.c.c."
\newline
Similarly (C3)(A3).  A natural approach is to consider $\langle N_s:s \in
[B]^{< \aleph_0} \rangle$ and use a subset $X \in [B]^{\text{otp}B}$ such
that for different uses we use almost disjoint $X$'s.  This was not completed
but we restrict ourselves to ``not only $P$ satisfies the $2^{\aleph_0}$-
c.c. but even $P^n$ (for each $n < \omega$).
\medskip

Concerning (C1) we cannot replace ${}^\mu n$ elements of
${\underset\tilde {}\to B}$ by 1, but we can use a directed system, so
``$P$ satisfies the $2^\mu$-c.c.", is replaced by
``for $\sigma < \mu,P^\sigma$ satisfies the $2^\mu$-c.c." (or slightly less).
\medskip

Another question is Velickovic's question answered for Borel c.c.c. forcing
in \cite{Sh:480}; i.e. (C4).
\bigskip

\centerline {Preliminaries}
\bigskip

\noindent
\subhead {0.A} \endsubhead  $Let  <^\ast_\chi$ be a well ordering of
\newline
$H(\chi ) = \{ x$:the transitive closure of $x$ has cardinality  $<\chi \}$  
agreeing with the usual well ordering of the ordinals, \newline
$P$  (and  $Q,R$)  
will denote forcing notions, i.e. partial order with a minimal element  
$\emptyset  = \emptyset_P$. 

A forcing notion  $P$  is $\lambda$-closed if every increasing sequence of 
members of  $P$,  of length less than  $\lambda$,  has an upper bound.
\bigskip

\subhead {0.B} \endsubhead   For sets of ordinals,  $A$  and  $B$,
define  $H^{OP}_{B,A}$ as the 
maximal order preserving bijection between initial segments of $A$ and $B$,  
i.e., it is the function with domain  $\{\alpha \in A:otp(\alpha \cap A) 
< otp(B)\}$  and  $H^{OP}_{A,B}(\alpha) = \beta$ if and only if  $\alpha  
\in  A$,  $\beta  \in  B$  and  $otp(\alpha \cap A) = otp(\beta \cap B)$.
\bigskip

\definition{Definition 0.1}  $\lambda \rightarrow^+ (\alpha)^{<\omega}_\mu$
holds provided that: \underbar{if} whenever  $F$  is a function from  
$[\lambda]^{<\omega}$ to $\lambda,F(w) < \text{ min}(w)$,
$C \subseteq  \lambda$  is a club \underbar{then} there is  
$A \subseteq  C$  of order type  $\alpha$  such that 
$\biggl[ w_1,w_2 \in  [A]^{<\omega},|w_1| =
|w_2| \Rightarrow  F(w_1) = F(w_2) \biggr]$. \newline
(See \cite{Sh:f},XVII,4.x).
\enddefinition
\bigskip

\remark{0.1A Remark}  1) If $\lambda$ is Ramsey cardinal then
$\lambda \rightarrow^+(\lambda)^{< \omega}_\mu$. \newline
2)  If $\lambda = \text{ Min}\{\lambda:\lambda \rightarrow (\alpha)^{< \omega}
_\mu\}$ then $\lambda$ is regular and $\lambda \rightarrow^+
(\alpha)^{< \omega}_\mu$.
\endremark

\definition{Definition 0.2}  $\lambda \rightarrow [\alpha]^n_{\kappa,\theta}$
if for every function  $F$  from  $[\lambda]^n$ to  $\kappa$  there is  $A 
\subseteq  \lambda $  of order type  $\alpha$  such that  $\{F(w):w \in  
[A]^n\}$  has power  $\le \theta$.
\enddefinition
\bigskip

\definition{Definition 0.3}  A forcing notion  $P$  satisfies the
Knaster condition (has 
property  $K$)  if for any  $\{p_i:i < \omega_1\} \subseteq P$  there is an 
uncountable  $A \subseteq \omega_1$ such that the conditions  $p_i$ and $p_j$ 
are compatible whenever  $i$, $j \in  A$.
\enddefinition
\newpage

\head {\S1 Consistency of ``c.c.c. $\times \, 2^{\aleph_0}$-c.c. =
$2^{\aleph_0}$ - c.c."} \endhead
\bigskip

The $a_i$'s are not really necessary but (hopefully) clarify.
\definition{1.1 Definition}  1)  ${\Cal K}_{\mu,\kappa}$ is the family
of  $\bar Q = \langle P_\gamma,{\underset\tilde {}\to Q_\beta},a_\beta:\gamma
\le \alpha,\beta < \alpha \rangle$, where:
\medskip
\roster
\item "{(a)}"  $\langle P_\gamma,{\underset\tilde {}\to Q_\beta}
:\gamma \le \alpha,
\beta < \alpha \rangle$  is a finite support iteration
\item "{(b)}"  every  $P_\gamma$, ${\underset\tilde {}\to Q_\gamma}$ satisfies
the c.c.c.
\item "{(c)}"  ${\underset\tilde {}\to Q_\beta}$ is a
$P_\beta$-name which depends
just on  $G_{P_\beta} \cap P^\ast_{a_\beta}$ (see below; hence it is in
$V[G_{P^*_\beta}])$, and  $|{\underset\tilde {}\to Q_\beta}| \le \kappa$
 and its set of members  $\subseteq V$  (for simplicity)
\item "{(d)}"  $a_\beta  \subseteq  \beta$,  $|a_\beta| \le \mu$ and $\gamma  
\in a_\beta \Rightarrow  a_\gamma \subseteq  a_\beta$.
\endroster
\medskip
\noindent
2)  For such  $\bar Q$  we call  $a \subseteq \ell g(\bar Q)$,$\bar Q$-closed 
if  $[\beta \in a \Rightarrow  a_\beta \subseteq a]$  and let
$$
\align
P^\ast_a = P^{\bar Q}_a =: \biggl\{ p \in P_\alpha:&\text{Dom}(p) \subseteq
a \text{ and for all } \beta \in \text{ Dom}(p):p(\beta) \in V \\
  &\text{(not a name) and } p \restriction a_\beta \Vdash ``p(\beta) \in 
Q_\beta" \biggr\}
\endalign
$$
\noindent
(so we are defining  $P^\ast_a$ by induction on  $\text{sup}(a)$) ordered
by the order of  $P_{\text{sup}(a)}$. \newline
3)  ${\Cal K}^k_{\mu,\kappa}$ is the class of
$\bar Q \in {\Cal K}_{\mu,\kappa}$ such that if
$\beta < \gamma \le \ell g(\bar Q)$,$cf(\beta) \ne \aleph_1$ then
$P_\gamma / P_\beta$ satisfies the Knaster condition (actually 
we can use somewhat less).  Let ${\Cal K}^n_{\mu,\kappa} = 
{\Cal K}_{\mu,\kappa}$. \newline
4)  If defining $\bar Q$ we omit $P_\alpha$ we mean $\dsize \bigcup_{\beta <
\alpha} P_\beta$ if $\alpha$ is limit, $P_\beta * {\underset\tilde {}\to
Q_\beta}$ if $\alpha = \beta + 1$. \newline
5)  We do not lose, if we assume ${\underset\tilde {}\to Q_\beta} \subseteq
[\kappa]^{< \aleph_0}$ and the order $\subseteq$; (then 1.2(1)(g) becomes 
trivial
as for closed $p,q \in P^*_j,p \restriction a \le q \restriction a$).
\enddefinition
\bigskip

\proclaim{1.2 Claim}  1) Assume  $x \in \{ n,k \}$ and
$\bar Q = \langle P_\gamma,{\underset\tilde {}\to Q_\beta},
a_\beta:\beta < \alpha,
\gamma \le \alpha \rangle \in  {\Cal K}^x_{\mu,\kappa}$.
\underbar{Then}
\medskip
\roster
\item "{(a)}"  for  $\alpha^\ast < \alpha$,  $\bar Q \restriction \alpha^\ast
 =: \langle P_\gamma,{\underset\tilde {}\to Q_\beta},
a_\beta:\beta < \alpha^\ast,
\gamma \le \alpha^\ast \rangle$  belongs to  ${\Cal K}^x_{\mu,\kappa}$ 
\item "{(b)}"  $P^\ast_\alpha$ is a dense subset of  $P_\alpha$ 
\item "{(c)}"  for any $\bar Q$-closed  $a \subseteq \alpha$, $P^\ast_a 
\lessdot P_\alpha$ (in particular  $P^\ast_\alpha$ is a dense subset 
of  $P_\alpha$);  moreover, if  $p \in P^\ast_\alpha$ then
$p \restriction a \in  P^\ast_a$ and \newline
$[p \restriction a \le q \in  P^\ast_a \Rightarrow  r =: q 
\cup  p \restriction (\alpha \backslash a) \in  P_\alpha \and p \le r
\and q \le r]$
\item "{(d)}"  for a $\bar Q$-closed  $a \subseteq \alpha$,  
$\langle P^\ast_{a \cap \gamma},{\underset\tilde {}\to Q_\beta},
a_\beta:\beta  \in  a,\gamma \in a \rangle$  
belongs to  ${\Cal K}^x_{\mu,\lambda}$ \newline
(except renaming; not used)
\item "{(e)}"  if  ${\underset\tilde {}\to Q_\alpha}$
is a $P^\ast_a$-name of a
c.c.c. forcing notion of cardinality $\le \kappa$, each member of
${\underset\tilde {}\to Q_\alpha}$ is from $V$, $a \subseteq \alpha$  
is $\bar Q$-closed,
$|a| \le \mu$  and  $P_{\alpha +1} = P_\alpha \ast 
{\underset\tilde {}\to Q_\alpha}$ and 
${\underset\tilde {}\to Q_\alpha}$ satisfies the Knaster condition or
at least $\beta < \alpha \Rightarrow
P_\alpha * {\underset\tilde {}\to Q_\alpha}/P_{\beta + 1}$ satisfies the
Knaster condition \underbar{then} \newline
$\langle P_\gamma,{\underset\tilde {}\to Q_\beta},
a_\beta:\beta < \alpha + 1,
\gamma \le \alpha  + 1 \rangle \in {\Cal K}^x_{\mu,\lambda}$
\item "{(f)}"  if $n < \omega,p_1,\dotsc,p_n \in P_{\alpha^*}$ and
{\roster
\itemitem{ $(*)$ }  for every $\beta \in \dsize \bigcup^n_{\ell = 1}
\text{ Dom}(p_\ell)$ for some $m = m_{\beta,\ell} \in \{ 1,\dotsc,n \}$
we have
$p_m \restriction \beta \Vdash$ ``$p_\ell(\beta) \le_{{\underset\tilde {}\to
Q_\beta}} p_m(\beta)$ for $\ell = \{1,\dotsc,n \}$"
\endroster}
\underbar{then} $p_1,\dotsc,p_n$ has a least common upper bound $p$ which is
defined by: \newline
Dom$(p) = \dsize \bigcup^n_{\ell = 1} \text{ Dom}(f),p_\ell(\beta)
= p_{m_{\beta,\ell}}(\beta)$, so in particular $p \in P_{\alpha^*}$ and
$\dsize \bigwedge^n_{\ell = 1} p_\ell \in P^*_{\alpha^*} \Rightarrow
p \in P^*_{\alpha^*}$
\item "{(g)}"  if $p_\ell \le p$ and
$p_\ell \in P^*_\gamma$ for $\ell < n$, and $a_k$ is $\bar Q$-closed for
$k < m$ \underbar{then} there is $p' \in P^*_\gamma$, such that $p \le p'$
and $P^*_{a_k} \models p_\ell \restriction a_k \le p' \restriction a_k$ for
$\ell < n,k < m$.
\endroster
\medskip

\noindent
2)  If  $x \in \{ n,k \}$ and $\delta < \lambda$  is a limit ordinal,
for  $\alpha < \delta$ we have \newline
$\langle P_\gamma,
{\underset\tilde {}\to Q_\beta},
a_\beta:\beta < \alpha,\gamma \le \alpha \rangle
 \in {\Cal K}^x_{\mu,\lambda}$ and $P_\delta = \dsize \bigcup_{\gamma <
\delta} P_\gamma$ then  $\langle P_\gamma,
{\underset\tilde {}\to Q_\beta},
a_\beta:\beta < \delta,\gamma \le \delta \rangle$
belongs to  ${\Cal K}^x_{\mu,\lambda}$.
\endproclaim
\bigskip

\demo{Proof}  Straightforward.
\enddemo
\bigskip

\noindent
Essentially by \cite{Sh:289},2.4(2),p.176 (which is slightly weaker and its
proof left to the reader, so we give details here).
\proclaim{1.3 Claim}  Assume  $\lambda \rightarrow^+ (\omega \alpha^*)
^{< \omega}_\mu$ (e.g. $\lambda$ a Ramsey cardinal, $\alpha^* = \lambda$)
$\chi > \lambda$, \newline
$x \in  H(\chi)$. \newline
1) There is an end extension strong $(\chi,\lambda,\alpha^*,\mu,\aleph_0,
\omega)$-system for $x$ (see Definition 1.3A). \newline
2)  There is an end extension $(\chi,\lambda,\alpha,\mu,\aleph_0,\omega)$-
system for $x$ \underbar{if} $x$ is Ramsey or $\lambda = \text{ Min}\{
\lambda:\lambda \rightarrow (\omega \alpha^*)^{< \omega}_\mu\}$ (also then
the condition holds for every $\mu' < \mu$).
\endproclaim
\bigskip

\definition{1.3A Definition}  1) We say $\bar N = \langle N_s:s \in [B]
^{< 1+n} \rangle$ is a $(\chi,\lambda,\alpha,\theta,\sigma,n)$-system if:
\medskip
\roster
\item "{(a)}"  $N_s \prec (H(\chi),\in)$ (or of some expansion) \newline
$\theta + 1 \subseteq N_s,\| N_s \| = \theta,{}^{\sigma >}(N_s) \subseteq
(N_s)$
\item "{(b)}"  $B \subseteq \lambda,\text{otp}(B) = \alpha$
\item "{(c)}"  $n \le \omega$ (equally is allowed but $1 + \omega = \omega$
so $s$ is always finite)
\item "{(d)}"  $N_s \cap N_t \subseteq N_{s \cap t}$
\item "{(e)}"  $N_s \cap B = s$
\item "{(f)}"  if  $|s| = |t|$  then  $N_s \cong  N_t$ say $H_{s,t}$ is an
isomorphism from  $N_s$ onto  $N_t$ (necessarily $H_{s,t}$ is unique)
\item "{(g)}"  \underbar{if}  $s' \subseteq  s$, $t' = \{\alpha \in t:
(\exists \beta  \in  s')[|\beta \cap  s| = |\alpha \cap t|]\}$
\underbar{then}  $H_{s',t'},H_{s,t}$ are compatible functions;
$H_{s,s} = id$, \newline
$H_{s,t} \supseteq H^{OP}_{s,t}$,
$H_{s_0,s_1} \circ  H_{s_1,s_2} = H_{s_0,s_2},H_{t,s} = (H_{s,t})^{-1}$
\item "{(h)}"  sup$(N_s \cap \lambda) < \text{ Min}\{ \alpha \in B:
\dsize \bigwedge_{\gamma \in s} \gamma < \alpha \}$.
\endroster
\medskip

\noindent
2)  We add the adjective ``strong" if in strengthen clause (d) by
\medskip
\roster
\item "{(d)$^+$}"  $N_s \cap N_t = N_{s \cap t}$ (so in clause (g),
$H_{s',t'} \subseteq H_{s,t}$).
\endroster
\medskip

\noindent
3)  We add the adjective ``end extension" if
\medskip
\roster
\item "{(i)}"  $s \triangleleft t \Rightarrow N_s \cap \lambda \triangleleft
N_t \cap \lambda$ (where $A \triangleleft B$) means 
$A = B \cap \text{ min}(B \backslash A)$
\endroster
\medskip

\noindent
4)  We add ``for $x$" if $x \in N_s$ for every $s \in [B]^{< 1+n}$, and
$H_{s,t}(x) = x$.
\enddefinition
\bigskip

\remark{1.3B Remark}  If $\lambda$ is a Ramsey cardinal (or much less see
\cite{Sh:f},XVII,4.x,\cite{Sh:289},\S4) then we have if 
$\gamma \in s \cap t,s \cap \gamma = t
\cap \gamma$ and $y \in N_s$ then in $(H(\chi),\in,<^*_\chi)$ the elements
$y$ and $H_{t,s}(y)$ realizes the same type over $\{ i:i < \gamma \}$.
[prove?]
\endremark
\bigskip

\demo{Proof}  1) Let $C = \{ \delta < \lambda:\text{for every } 
\alpha < \delta
\text{ there is } N \prec (H(\chi),\in,<^*_\chi)$ such that $\mu + 1 + 
\alpha \subseteq
N \text{ and } \text{ sup}(N \cap \lambda) < \delta \}$. Clearly $C$ is a
club of $\lambda$. \newline

Let  $B_0 = \{\alpha_i:i < \omega \alpha^* \} \subseteq C$,  
$(\alpha_i$ strictly increasing) be indiscernible in  \newline
$(H(\chi),\in,<^\ast_\chi,x)$ (see Definition 0.1).  
Let  $B = \{\alpha_i:i < \omega \alpha^* \text{ limit}\}$.  
For  $s \in [B_0]^{<\aleph_0}$ let $N^0_s$ = the Skolem 
Hull of $s \cup \{ i:i \le \mu \} \cup \{ x,\lambda \}$
under the definable functions of $(H(\chi),\in,<^*_\chi)$ and \newline
$N_s = \cup \biggl\{ N^0_{t_1} \cap  N^0_{t_2}:t_1,t_2 \in 
[\{\alpha_i:i < \omega \alpha^* \}]^{<\aleph_0} \text{ and }
s = t_1 \cap t_2 \biggr\}$. \newline
Clearly
\medskip
\roster
\item "{$(*)$}"  $\| N_s \| \le \mu$ and $\{ x,\lambda \} \subseteq N_s$.
\endroster
\medskip

\noindent
Now we shall show
\medskip
\roster
\item "{$(*)_1$}"  if $s \in [B]^{< \aleph_0},y \in N_s$ \underbar{then}
for every finite $t \subseteq B_0$ there is $u \in [B_0]^{<\aleph_0}$ such
that $s \subseteq u,u^* \cap t \subseteq s$ and $y \in N^0_u$.
\endroster
\medskip

\noindent
As $y \in N_s$ there are $s_1,s_2 \in [B_0]^{< \aleph_0}$ such that
$y \in N^0_{s_1} \cap N^0_{s_2}$ and $s=s_1 \cap s_2$.  Let $s_1 \cup s_2 =
\{ \alpha_{i_0},\dotsc,\alpha_{i_{m-1}} \}$ (increasing), and let \newline
$n^* = 
\text{ sup}\{n:
\text{for some } \beta,\beta + n \in t \} + 1$, and define for $\ell \le m$ a
function $f_\ell$ with domain $s_1 \cup s_2$, such that
\enddemo

$$
f_\ell(\alpha_{i_k}) = \left\{
\alignedat2
&\alpha_{i_k+1} \qquad &&\text{\underbar{if}} 
\quad k \ge m - \ell \text{ and } i_k \notin s \\
&\alpha_{i_k}    &&\text{\underbar{otherwise}}
\endalignedat
\right.
$$

\noindent
Note that
\medskip
\roster
\item "{$\bigotimes_1$}"  for $\ell < m,f_\ell \restriction s_1 = 
f_{\ell + 1} \restriction s_1$ or $f_\ell \restriction s_2 = f_{\ell + 1}
\restriction s_2$ (or both) \newline
[why? as $i_\ell \in s_2 \backslash s_1 \backslash s$ or
$i_\ell \in s_2 \backslash s_1 \backslash s$
or $i_\ell \in s = s_1 \cap s_2$].
\item "{$\bigotimes_2$}"  $f_\ell$ is order preserving with domain
$s_0 \cup s_1,f_\ell \restriction s =$ the identity.
\endroster
\medskip

\noindent
As $y \in N^0_{s_1} \cap N^0_{s_2}$ there are terms $\tau_1,\tau_2$ such that
$$
y = \tau_1(\dotsc,\alpha_{i_k},\ldots)_{\alpha_{i_k} \in s_1} = 
\tau_2(\dotsc,\alpha_{i_k},\ldots)_{\alpha_{i_k} \in s_2}.
$$

\noindent
Using the indiscernibility of $B_0$ we can prove by induction on $\ell \le m$
that
\medskip
\roster
\item "{$\bigotimes_{3,\ell}$}"  $y = \tau_1(\dotsc,f_\ell(i_{\alpha_{i_k}},
\ldots)_{\alpha_{i_k} \in s_1} = \tau_2(\dotsc,f_\ell(\alpha_{i_k}),
\ldots)_{\alpha_{i_k} \in s_2}$.
\endroster
\medskip

\noindent
[Why?  For $\ell = 0$ this is given by the choice of $\tau_1,\tau_2$.
For $\ell + 1$ note that by $\otimes_2$, $f_{\ell + 1} \circ f^{-1}_\ell$ is
an order preserving function from $\text{Rang}(f_\ell)$ onto
$\text{Rang}(f_{\ell + 1}$). \newline
By $\otimes_{3,\ell}$ and ``$B_0$ is indiscernible" we know \newline
$\tau_1(\dotsc,f_\ell(\alpha_{i_k}),\ldots)
_{\alpha_{i_k} \in s_1} = \tau_2(\dotsc,f_\ell(\alpha_{i_k}),\ldots)
_{\alpha_{i_k} \in s_2}$.  By the
last two sentences and the indiscernibility of $B_0$

$$
\tau_1(\dotsc,(f_{\ell + 1} \circ f^{-1}_\ell)(f_\ell(\alpha_{i_k})),\ldots)
_{\alpha_{i_k} \in s_1} = \tau_2(\dotsc,(f_{\ell + 1} \circ f^{-1}_\ell)
(f_\ell(\alpha_{i_k})),\ldots)_{\alpha_{i_k} \in s_2}.
$$

\noindent
But $(f_{\ell + 1} \circ f^{-1}_\ell)(f_\ell(\alpha_{i_k})) = f_{\ell + 1}
(\alpha_{i_k})$ so

$$
\tau_1(\dotsc,f_{\ell + 1}(\alpha_{i_k}),\ldots)_{\alpha_{i_k} \in s_1} = 
\tau_2(\dotsc,f_{\ell + 1}(\alpha_{i_k}),\ldots)_{\alpha_{i_k} \in s_2}.
$$

\noindent
But by $\otimes_1$ for some $e \in \{ 1,2 \}$ we have 
$f_\ell \restriction s_e =
f_{\ell + 1} \restriction s_e$, so \newline
$\tau_e(\dotsc,f_{\ell + 1}(\alpha_{i_k}),\ldots)_{\alpha_{i_k} \in s_e} = 
\tau_e(\dotsc,f_\ell(\alpha_{i_k}),\ldots)_{\alpha_{i_k} \in s_e}$ 
but the latter is equal to
$y$ (by the induction hypothesis), hence the former so by the last sentence

$$
y = \tau_1(\dotsc,f_{\ell + 1}(\alpha_{i_k}),\ldots)_{\alpha_{i_k} \in s_1} = 
\tau_2(\dotsc,f_{\ell + 1}(\alpha_{i_k}),\ldots)_{\alpha_{i_k} \in s_2}.
$$

\noindent
So we have caried the induction on $\ell \le m$, and for $\ell = m$ we get
$y \in N^0_{f_m(s_1)}$, but by the choice of $n^*$ and $f_m$ clearly
$f_m(s_1) \cap t \subseteq s$,  and we have proved $(*)_1$. \newline
Now we can note
\medskip
\roster
\item "{$(*)_2$}"  if $s \in [B]^{< \aleph_0}$ and 
$y_1,\dotsc,y_n \in N_s$ then for some
$s_1,s_2 \in [B_0]^{< \aleph_0}$ we have: $s = s_1 \cap s_2$ and $y_1,
\dotsc,y_n \in N^0_{s_1} \cap N^0_{s_2}$.
\endroster
\medskip

\noindent
[Why?  We can find $u_1,\dotsc,u_n \in [B_0]^{<\aleph_0}$ such that
$s \subseteq u_\ell,
y_\ell \in N^0_{u_\ell}$ (as $y_\ell \in N_s$).  Now by $(*)_1$ for each
$\ell = 1,2,\dotsc,n$ we can find $v_\ell \in [B_0]^{< \aleph_0}$ such that
$s \subseteq v_\ell, s=v_\ell \cap (\dsize \bigcup^n_{m=1} u_m)$ and
$y_\ell \in N^0_{v_\ell}$.  Let $u = \dsize \bigcup^n_{i=1} u_\ell,
v = \dsize \bigcup^n_{\ell = 1} u_\ell$, clearly
$y_1,\dotsc,y_n \in N^0_u \cap N^0_v$ and $u \cap v = s$, as required].
\newline
Now as we have Skolem functions $(*)_2$ implies
\medskip
\roster
\item "{$(*)_3$}"  $N_s \prec (H(\chi),\in,<^*_\chi)$
\endroster
\medskip

\noindent
Also trivially
\medskip
\roster
\item "{$(*)_4$}"  $N^0_s \prec N_s$ hence $\mu + 1 \subseteq N_s$
\medskip
\item "{$(*)_5$}"  $s \subseteq t \Rightarrow N_s \prec N_t$.
\endroster
\medskip

\noindent
Also
\medskip
\roster
\item "{$(*)_6$}"  $N_{s_1} \cap N_{s_2} = N_{s_1 \cap s_2}$ for
$s_1,s_2 \in [B]^{< \aleph_0}$.
\endroster
\medskip

\noindent
[Why?  The inclusions $N_{s_1 \cap s_2} \subseteq N_{s_1} \cap N_{s_2}$
follows from $(*)_5$; for the other direction let $y \in N_{s_1} \cap
N_{s_2}$.  By $(*)_1$ as $y \in N_{s_1}$ there is $t_1$ such that
$s_1 \subseteq t_1 \in [B_0]^{< \aleph_0},t_1 \cap (s_1 \cup s_2) = s_2$
and $y \in N^0_{t_1}$.  By $(*)_1$, as $y \in N_{s_2}$ there is $t_2$
such that $s_2 \subseteq t_2 \in [B_0]^{< \aleph_0},
t_2 \cap (s_1 \cup s_2 \cup t_1) = s_1$ and $y \in N^0_{t_2}$.  So
$y \in N^0_{t_1} \cap N^0_{t_2}$, but easily
$t_1 \cap t_2 = s_1 \cap s_2$].
\medskip
\roster
\item "{$(*)_7$}"  $\text{sup}(N_s \cap \lambda) < \text{ Min}
\{\alpha \in B:\dsize \bigwedge_{\gamma \in s} \gamma < \alpha\}$. \newline
[why?  as $B_0 \subseteq C$ and see the Definition of $C$].
\endroster
\medskip

\noindent
Now check that (a)-(h) of Definition 1.3A holds. \newline
Now $\langle N_s:s \in [B]^{<\aleph_0} \rangle$ is as required. \newline
2)  If $\lambda$ is Ramsey, without loss of generality $\text{otp}(B_0)
= \lambda$ and it is easy to check 1.3A(i).  The other case is like
\cite{Sh:289},\S4. \hfill$\square_{1.3}$
%\enddemo
\bigskip

\proclaim{1.4 Theorem}  Assume  $\aleph_0 < \mu \le \kappa < \lambda =
cf(\lambda)$,$\lambda$  strongly inaccessible,  $\lambda$  a Ramsey cardinal,
and  $\diamondsuit_{\{\delta <\lambda:cf(\delta) = \aleph_1\}}$ (can be
added by a preliminary forcing). \newline
\underbar{Then} we have  $P$  such that:
\medskip
\roster
\item "{(a)}"  $P$  is a c.c.c. forcing of cardinality $\lambda$ adding
$\lambda$  
reals (so the cardinals and cardinal arithmetic in  $V^P$ should be clear),
in 
particular  in $V^P$ we have $2^{\aleph_0} = \lambda$
\item "{(b)}"  $\Vdash_P$ ``$MA$  holds for c.c.c. forcing notions
of cardinality  $\le \mu$ and  $< \lambda$  dense sets
(and even for c.c.c. forcing notions of
cardinality  $\le \kappa$  which are from  $V[A]$  for some
$A \subseteq  \mu$)"
\item "{(c)}"  $\Vdash_P$ ``if $B$ is a $\lambda$-c.c. Boolean algebra,
$x_i \in  B \backslash \{0\}$  for  $i < \lambda$  \underbar{then} for
some  $Z \subseteq \lambda $,  $|Z| = \aleph_1$ and  $\{x_i:i \in  Z\}$
generates a proper filter of  $B$  (i.e. no finite intersection is $0_B$)"
\item "{(d)}"  $\Vdash_P$ ``if $B_1$ is a c.c.c. Boolean algebra, $B_2$ is a
$\lambda$-c.c. Boolean algebra then  $B_1 \times B_2$ is a
$\lambda$-c.c. Boolean algebra."
\endroster
\endproclaim
\bigskip

\demo{Proof}  Let  $\langle A_\delta:\delta < \lambda,cf(\delta) = 
\aleph_1 \rangle$  exemplifies the diamond.  We choose by induction on  
$\alpha  < \lambda$,  $\bar Q^\alpha = \langle P_\gamma,
{\underset\tilde {}\to Q_\beta},
a_\beta:\gamma \le \alpha,\beta < \alpha \rangle 
\in {\Cal K}^n_{\mu,\kappa}$ such that $\alpha^1 < \alpha \Rightarrow  
\bar Q^{\alpha^1} = \bar Q^\alpha \restriction \alpha^1$.  In limits
$\alpha$  use 1.2(2), for  $\alpha = \beta + 1$, $cf(\beta) \ne \aleph_1$
take care of (b) by suitable bookkeepping using 1.2(1)(e). 
If  $\alpha = \beta + 1$,
$cf(\beta) = \aleph_1$ and  $A_\beta$ codes  $p \in P_\beta$ and
$P_\beta $-names of a Boolean algebra
${\underset\tilde {}\to B_\beta}$
and sequence  $\langle {\underset\tilde {}\to x^\beta_i}:
i < \beta \rangle$  of non-zero members of ${\underset\tilde {}\to B_\beta}$,
and 
$p$ forces ($\Vdash_{P_\beta}$)  that there is in
$V[{\underset\tilde {}\to G_{P_\beta}}]$ some c.c.c. forcing notion  $Q$  of 
cardinality  $\le \mu$  adding some  $Z \subseteq \beta$,  $|Z| = 
\aleph_1$ with  $\{x^\beta_i:i \in  Z\}$  generating a proper filter
of  ${\underset\tilde {}\to B_\beta}$ \underbar{then} we choose
${\underset\tilde {}\to Q_\beta}$,  if
  $p \in {\underset\tilde {}\to G_{p_\beta}}$,
as such  $Q$.  If  $p \notin {\underset\tilde {}\to G_{P_\beta}}$
 or there is no
such  ${\underset\tilde {}\to Q}$ in $V[{\underset\tilde {}\to G_{p_\beta}}]$,
\underbar{then}  ${\underset\tilde {}\to Q_\beta}$ is e.g. Cohen forcing. 

So every  $\bar Q^\alpha$ is defined, let
$P = \dsize \bigcup_{\gamma < \lambda} P_\gamma$.  Clearly ($\alpha$) + (b)
holds and (d) follows by (c).  So the rest of the proof is dedicated to
proving (c). \newline
\noindent
So let  $p \in  P$,  $p \Vdash ``{\underset\tilde {}\to B}$  a $\lambda$-c.c.
Boolean algebra,  ${\underset\tilde {}\to x_i} \in B \backslash \{0_B\}$
for $i < \lambda"$
without loss of generality the set of members of
${\underset\tilde {}\to B}$ is $\lambda$. 

Let  $x = \left< P,p,{\underset\tilde {}\to B},
\langle {\underset\tilde {}\to x_i}:i <
\lambda \rangle \right>$,  $\chi  = \lambda^+$,  by Claim 1.3 
there are  $A \in [\lambda]^\lambda$ and  $\langle N_s:s \in [A]^{<\aleph_0}
\rangle$  as there (for $\kappa = \mu + \kappa$ here standing for $\mu$
there).  Let
$$
\align
C = \biggl\{ \delta < \lambda:&\delta \text{ a strong limit cardinal }
> \kappa + \mu, [\alpha < \delta \Rightarrow \bar Q \restriction \alpha  \in 
H(\delta)], \\
 &\delta = \text{ sup}(A \cap \delta),s \in [A \cap \delta]^{<\aleph_0} 
\Rightarrow  \text{ sup}(\lambda \cap N_s) < \delta, \\
  &{\underset\tilde {}\to B} \restriction \delta \text{ a }
 P_\delta \text{-name, and for } i < \delta \text{ we have }
{\underset\tilde {}\to x_i} \text{ a } P_\delta \text{-name}
\biggr\}.
\endalign
$$

\noindent
For some  accumulation point  $\delta$  of $C$,  $cf(\delta) = \aleph_1$ and 
$A_\delta$ codes $\left< p,{\underset\tilde {}\to B} \restriction \delta,
\langle {\underset\tilde {}\to x_i}:i < \delta \rangle \right>$.  We shall
show that for some  $q$,  
$p \le q \in  P_\delta$ and  $q \Vdash_{P_\delta}$ ``there is  $Q$  as 
required above".  By the inductive choice of 
 ${\underset\tilde {}\to Q_\delta}$
this suffices. \newline
\noindent
Let  $A^\ast \subseteq  A \cap \delta$,  otp$(A^\ast) = \omega_1$,  
$\delta = \text{ sup}(A^\ast)$ and  $\langle \delta_i:i < \omega_1 \rangle$  
increasing continuous,  $\delta = \dsize \bigcup_{i < \omega_1} \delta_i$,
$\delta_i \in C$,  $A^\ast \cap \delta_0 = \emptyset$,
$|A^\ast \cap [\delta_i,\delta_{i+1})| = 1$. \newline
\noindent
In $V^{P_\delta}$ we define:

$$
{\underset\tilde {}\to Q}
 = \biggl\{ u:u \in [A^\ast]^{<\aleph_0}, \text{ and }  
{\underset\tilde {}\to B}
 \models ``\dsize \bigcap_{i \in u} {\underset\tilde {}\to x_i}
\ne {0_{\underset\tilde {}\to B}}" \biggr\}
$$
\medskip

\noindent
ordered by inclusion.  It suffices to prove that some  
$q$,  $p \le q \in P_\delta$,
$q$  forces that: ${\underset\tilde {}\to Q}$  is c.c.c. with
$\cup \, {\underset\tilde {}\to G_{\underset\tilde {}\to Q}}$
 an uncountable set;
now clearly $q$ forces that $\{ {\underset\tilde {}\to x_i}:i \in \cup  
{\underset\tilde {}\to G_{\underset\tilde {}\to Q}} \}$  generates a proper
filter of  ${\underset\tilde {}\to B}$.\newline
If not, we can find  $q_i$, $u_i$ such that:

$$
p \le q_i \in  P^*_\delta \text{ and } q_i \Vdash_{P_\delta} ``u_i \in  
{\underset\tilde {}\to Q}" \text{ (where }  u_i \in [A^\ast]^{<\aleph_0})
$$
\smallskip

\noindent
and  $\langle (q_i,u_i):i < \omega_1 \rangle$  are pairwise incompatible in  
$P_\delta \ast {\underset\tilde {}\to Q}$. \newline
\noindent
Let  $v_i$ be a finite subset of $A^\ast$ such that:  $u_i \subseteq v_i$,  
and
\medskip
\roster
\item "{$(*)$}"  $[v \subseteq  A^\ast \and v \text{ finite } \and
\gamma  \in (\text{Dom } q_i) \cap  N_v \Rightarrow \gamma \in (\text{Dom }
q_i) \cap N_{v \cap v_i}]$.
\endroster
\medskip 
\noindent
By Fodor's Lemma for some stationary, $S \subseteq \omega_1,u^*,v^*,n^*$
and  $i(\ast)$ we have: for  $i < j$  in  $S$,

$$
v_i \cap \delta_i = v^\ast \subseteq \delta_{i(*)}, v_i \subseteq
\delta_j,|v_i| = n^\ast,u_i \cap \delta_i = u^\ast 
$$

$$
i(*) = \text{ Min}(S)
$$

$$
\{|\gamma \cap v_i)|:\gamma \in u_i\} \text{ does not depend on } i
$$
 
$$
q_i \restriction \delta_i \in P^*_{\delta_{i(\ast )}}
$$

$$
q_i \in  P^*_{\delta_j}.
$$
\medskip
\noindent
Let  $b_i =: N_{v_i} \cap \lambda$, so  $b_i$ is necessarily 
$\bar Q^\delta $-closed and  $|b_i| = \kappa$.  Let
$q^1_i = q_i \restriction b_i$,  so necessarily
$q^1_i \in  P^\ast_{b_i}$ (see 2.2(1)(c)).  Easily  
$P^\ast_{b_i} \subseteq N_{v_i}$ (though do not belong to it) so
$q^1_i \in  N_{v_i}$. 
\medskip

Let  $q^2_i =: H_{v_{i(*)},v_i}(q^1_i)$,  so  $q^2_i
\in  P^\ast_{b_{i(*)}}$; let $q^3_i =: (q_i \restriction \delta_{i(*)}) \cup
\biggl[ q^2_i \restriction (N_{v_{i(*)}} \cap \lambda \backslash 
\delta_{i(*)}) \biggr]$ by 1.2(1)(c) we know $q^3_i \in P^*_{\text{sup}(b
_{i(*)})+1}$ and $q^2_i \le q^3_i$, even without loss of generality
$q^2_i \le q^3_i \restriction b_{i(*)}$.  As  $P^\ast_{\text{sup}b_{i(*)}+1}
\lessdot P_\delta$ and $P_\delta$ satisfies the c.c.c. clearly 
for some  $i < j$
from  $S$, $q^3_i,q^3_j$,  are compatible in
$P^\ast_{\text{sup}(b_{i(*)})+1}$,  so let  $r \in
P^\ast_{\text{sup}(b_{i(*)})+1}$ be a 
common upper bound.  So $q^3_i \restriction (\delta_{i(*)} \cap 
b_{i(*)}) \le r \restriction
(\delta_{i(*)} \cap b_{i(*)})$ and $q^3_j \restriction (\delta_{i(*)}
\cap b_{i(*)}) \le r \restriction (\delta_{i(*)} \cap b_{i(*)})$ and
$q^3_i \restriction b_{i(*)} \le r \restriction b_{i(*)},q^3_j \restriction
b_{i(*)}$. \newline
\medskip

\noindent
Without loss of generality $\text{Dom}(r) \subseteq b_{i(*)} \cup
\delta_{i(*)}$ (allowed as $b_{i(*)}$ and $\delta_{i(*)}$ are closed, see
1.2(1)(c)); let $r_i = H_{v_i,v_{i(*)}}(r \restriction b_{i(*)}$)
and similarly
$r_j = H_{v_j,v_{i(*)}}(r \restriction b_{i(*)})$.
\enddemo
\bigskip

\noindent
Note that $r_i \in P^*_{\delta_j},r_j \in P^*_\delta,r_j \restriction
\delta_j = r_i \restriction \delta_i = r \restriction \delta_{i(*)}$.  Hence
$r_i \cup r_j \in P^*_\delta$.
\subhead {Case 1} \endsubhead  $r_i \cup r_j$  do not force (i.e. $\Vdash
_{P_\delta}$) that

$$
{\underset\tilde {}\to B} \models ``\dsize \bigcap_{\alpha \in u_i \cup u_j}
{\underset\tilde {}\to x_\alpha} = {0_{\underset\tilde {}\to B}}".
$$

\medskip
\noindent
Then there is  $r' \in P_\delta$,  $r_i \le r'$, $r_j \le r'$  forcing
the negation.  
So without loss of generality 
$r' \in  P^\ast_\delta$, and (as all parameters appearing in the
requirements on $r'$ are in $N_{v_i \cup v_j}$ also)
$r' \in P^\ast_{\lambda \cap (N_{v_i\cup v_j})}$.  Now  \newline
$r',r,q_i,q_j$ has an upper bound  $r'' \in P_\delta$. \newline
[Why?  By 1.2(1)(f), we have to check the condition $(*)$ there, so let
\newline
$\beta \in \text{ Dom}(r') \cup \text{ Dom}(r) \cup \text{ Dom}(q_i)
\cup \text{ Dom}(q_j)$].
\bigskip

\subhead{Subcase a} \endsubhead  $\beta \in \delta_{i(*)} \backslash
N_{v_i \cup v_j}$.  Note that $N_{v_i \cup v_j} \cap \delta_{i(*)} = N_{v^*}
\cap \lambda = b_{i(*)}$ (see choice of the $N_u$'s and definition of the
$b_\varepsilon$'s) but $\text{Dom}(r') \subseteq N_{v_i \cup v_j} \cap
\lambda$, so $\beta \notin \text{ Dom}(r')$.  Now

$$
q_i \restriction \delta_i = q_i \restriction
\delta_{i(*)} = q^3_i \restriction \delta_{i(*)} \le r
$$

$$
q_j \restriction \delta_j = q_j \restriction
\delta_{i(*)} = q^3_j \restriction \delta_{i(*)} \le r.
$$
\medskip

\noindent
So $r \restriction \beta \Vdash_{P_\beta}$ ``$q_i(\beta) \le r(\beta),
q_j(\beta) \le r(\beta)$" and $\beta \notin \text{ Dom}(r')$.  So we have
confirmed $(*)$ from 1.2(1)(f) for this subcase.
\bigskip

\subhead{Subcase b} \endsubhead  $\beta \in \delta_{i(*)} \cap 
N_{v_i \cup v_j}$. \newline

Exactly as above: \newline
$N_{v_i \cup v_j} \cap \delta_{i(*)} = N_{v^*} \cap \lambda =
b_{i(*)}$, so $\beta \in N_{v^*},\beta \in \delta_{i(*)} \cap b_{i(*)}$.
Also \newline
$q_i \restriction b_{i(*)} = q^1_i \restriction \delta_{i(*)} =
q^2_i \restriction \delta_{i(*)} = q^3_i \restriction (\delta_{i(*)} \cap
b_{i(*)}) \le r \restriction (\delta_{i(*)} \cap b_{i(*)})$ and \newline
$q_j \restriction b_{i(*)} = q^1_j \restriction \delta_{i(*)} =
q^2_j \restriction \delta_{i(*)} = q^3_j \restriction (\delta_{i(*)} \cap
b_{i(*)}) \le r \restriction (\delta_{i(*)} \cap b_{i(*)})$ and \newline
$r \restriction (\delta_{i(*)} \cap b_{i(*)}) \le r'$ (as $H_{v_i,v_{i(*)}}$
is the identity on $\delta_{i(*)} \cap b_{i(*)}$).

The last three inequalities confirm the requirement in 1.2(1)(f) (as
$\beta \in \delta_{i(*)} \cap b_{i(*)}$, see above).
\bigskip

\subhead {Subcase c} \endsubhead  $\beta \in (\delta \backslash
\delta_{i(*)}) \backslash N_{v_i \cup v_j}$. \newline

In this case $\beta \notin \text{ Dom}(r')$ (as $r' \in N_{v_i \cup v_j}$).
Also $\delta_{i(*)} < \delta_i < \delta_j < \delta$ and:

$$
\text{Dom}(r) \backslash \delta_{i(*)} \subseteq (b_{i(*)} \cup \delta
_{i(*)}) \backslash \delta_{i(*)} \subseteq [\delta_{i(*)},\delta_i)
$$

$$
\text{Dom}(q_i) \backslash \delta_{i(*)} \subseteq [\delta_i,\delta_j)
$$

$$
\text{Dom}(q_j) \backslash \delta_{i(*)} \subseteq [\delta_j,\delta).
$$

So $\beta$ belongs to at most one of $\text{Dom}(r'),\text{Dom}(r),
\text{Dom}(q_i),\text{Dom}(q_j)$ so the requirement $(*)$ from 1.2(1)(f)
holds trivially.
\bigskip

\subhead{Subcase d} \endsubhead
$\beta \in (\delta \backslash \delta_{i(*)}) \cap
N_{v_i \cup v_j}$. \newline

Clearly $\beta \notin \text{ Dom}(r)$. \newline

We know $q_i \restriction b_i = q^1_i,r_i \le r',
H_{v_{i(*)},v_i}(q^1_i) = q^2_i \le q^3_i \restriction b_{i(*)} \le r
\restriction b_{i(*)}$ hence $q^1_i \le H^{-1}_{v_{i(*)},v_i}(r \restriction
b_{i(*)}) = H_{v_i,v_{i(*)}}(r \restriction b_{i(*)}) = r_i$
but $r_i \le r'$, so together $q^1_i \le r'$, and similarly $q^1_j \le r'$.
As we have noted $\beta \notin \text{ Dom}(r)$ we have finished confirming
condition $(*)$ from 1.2(1)(f).

So really $r,r',q_i,q_j$ has a least common upper bound,  
hence  $(r'',u_i \cup u_j) \in P_\delta \ast {\underset\tilde {}\to Q}$
exemplified  $(q_i,u_i)$,  $(q_j,u_j)$  are compatible, as required.
\bigskip

\subhead {Case 2} \endsubhead  Not 1. 

Let  $\langle s_\beta :\beta  < \lambda \rangle $  be such that:

$$
s_\beta  \in  [A]^{<\aleph_0}, v^\ast \subseteq s_\beta,  
|s_\beta \backslash v^\ast| = |v_i \backslash v^\ast|,\text{sup}(v^*) <
\delta_{i(*)} < \text{ min}(s_\beta \backslash v^*)
$$

$$
\delta < \text{ min}(s_\beta \backslash v^*) \text{ (for simplicity)}
$$
 
$$
\beta < \gamma \Rightarrow  \text{ max}(s_\beta) < 
\text{ min}(s_\gamma \backslash v^\ast).
$$
\medskip

\noindent
As the truth value of $\dsize \bigcap_{\alpha \in u_i} x_\alpha$ is a
$P^*_a$-name for some closed $a \in N_{v_i}$ of cardinality $\le \mu$,
and $q_i \Vdash [{\underset\tilde {}\to B} \models$ ``$\dsize
\bigcap_{\alpha \in u_i} x_\alpha \ne 0_{\underset\tilde {}\to B}"]$ clearly
$q^1_i \Vdash [{\underset\tilde {}\to B} \models$ ``$\dsize
\bigcap_{\alpha \in u_i} x_\alpha \ne 0_{\underset\tilde {}\to B}"]$.
\medskip

\noindent
For $\beta < \lambda$ let
$\gamma^\beta = H_{s_\beta,v_{i(*)}}(r \restriction b_{i(*)}$,
and $u'_\beta = H_{s_\beta,v_{i(*)}}(u_0)$.  Let

$$
{\underset\tilde {}\to Y} = \{\beta < \lambda:r^\beta  \in  
{\underset\tilde {}\to G_P} \}.
$$

\medskip
\noindent
Clearly:

$$
r^\beta \Vdash_{P_\lambda}[{\underset\tilde {}\to B} \Vdash
``\dsize \bigcap_{i \in u'_\beta} {\underset\tilde {}\to x_i} \ne
{0_{\underset\tilde {}\to B}}"].
$$

\medskip
\noindent
Clearly  $p \le r^\beta$ and for some  $\beta$ we have $r^\beta \Vdash
``{\underset\tilde {}\to Y} \in  [\lambda]^\lambda$
 (and  $p \in {\underset\tilde {}\to G_P}$)"  and by the assumption of the
case: 

$$
\align
p \Vdash &``\bigg\{ \dsize \bigcap_{i \in u'_\beta} x_i:\beta \in Y \biggr\}
\text{ is a set of non-zero members of } {\underset\tilde {}\to B} \\
  &\text{ any two having zero intersection in } {\underset\tilde {}\to B}".
\endalign
$$

\medskip
\noindent
This contradicts an assumption on  $B$. \hfill$\square_{1.4}$
\bigskip

\centerline {$* \qquad \qquad * \qquad \qquad *$}
\bigskip
 
\noindent
We can phrase the consistency result as one on colouring.
\bigskip

\proclaim{1.5 Lemma}  1) In 1.4 we can add: 
\medskip
\roster
\item "{(e)}"  if  $c$  is a symmetric function from  
$\biggl[ 2^{\aleph_0} \biggr]^{<\omega}$ to $\{ 0,1 \}$ then
at least one of the following holds:   
{\roster
\itemitem{ $(\alpha)$ }  we can find pairwise disjoint
$w_i \subseteq 2^{\aleph_0}$ for  $i < 2^{\aleph_0}$ such that: \newline
$c \restriction [w_i] < \aleph_0$ is constantly zero but  \newline
$\dsize \bigwedge_{i < j} (\exists u \subseteq w_i,\exists v \subseteq w_j)
\biggl[ c[u \cup v] = 1 \biggr]$
\itemitem{ $(\beta)$ }  we can find an unbounded $B \subseteq 2^{\aleph_0}$
such that  $c \restriction [B]^{<\omega}$ is constantly $0$.
\endroster}
\endroster
\endproclaim
\bigskip

\noindent
It is natural to ask:

\subhead {1.6 Question} \endsubhead  Can we replace  $2^{\aleph_0}$ by
$\lambda < 2^{\aleph_0}$? \, $\aleph_1$ by  $\mu < \lambda$?  What is the
consistency strength of the statements we prove consistent? (see later).
Does $\lambda$ strongly inaccessible $k^2_2$-Mahlo (see [Sh289]) suffice?
\bigskip

\demo{1.7 Discussion}  Of course, 1.5(e) $\Rightarrow$ 1.4(c) $\Rightarrow$
1.4(d).  Starting with $\lambda$ weakly compact we can get a c.c.c. forcing
notion  $P$ of cardinality  $\lambda$,
such that in  $V^P,2^{\aleph_0} = \lambda$ and (e) of 1.5 holds for
$c:\biggl[ 2^{\aleph_0} \biggr]^2 \rightarrow  \{0,1\}$  (so 
$c(u) = 0$  if  $|u| \ne 2)$  and this suffices 
for the result.  Also we can generalize to higher cardinals.  We shall
deal with this elsewhere. 
\enddemo
\bigskip

\proclaim{1.8 Theorem}  Concerning the consistency strength, in 1.4 it
suffices to assume
\medskip
\roster
\item "{$(*)$}"  $\lambda$ is strongly inaccessible and for every
$F:[\lambda]^{< \aleph_0} \rightarrow \mu$ and club $C$ we can find
$B \subseteq C$, (or just $B \subseteq \lambda$) \, $\text{otp}(B) =
\omega_1$ such that
{\roster
\itemitem{ (a) }  $B$ is $F$-indiscernible i.e. if $n < \omega,u,v \in
[B]^n$ then $F(u) = F(v)$
\itemitem{ (b) }  for every $n < \omega$ there is $B' \in [C]^\lambda$
such that:
$$
\text{if }u \in [B']^n \text{ and } v \in [B]^n \text{ then }
F(u) = F(v)
$$
\endroster}
\endroster
\endproclaim
\bigskip

\demo{Proof}  Let $R = \{ \bar Q:\bar Q \in H(\lambda),\bar Q \in {\Cal K}
^n_{\mu,\kappa}\}$ ordered by $\bar Q^1 < \bar Q^2$ if \newline
$\bar Q^1 = \bar Q^2 \restriction \ell g(\bar Q^1)$.  Clause (b) takes care
also of ``the end extension" clause and for 1.3(A)(4), Clause (b) the proof
is the same. \newline
A somewhat less natural property though suffices. \newline
(Note: Clause (b) also helps to get rid of the club $C$).
\enddemo
\bigskip

\proclaim{1.9 Claim}  In 1.4 it suffices to assume
\medskip
\roster
\item "{$(*)'$}"  if $F:[\lambda]^{< \aleph_0} \rightarrow \mu$ then there
is $B \subseteq \lambda,\text{otp}(B) = \omega$ such that
{\roster
\itemitem{ (a) }  $F \restriction [B]^n$ is constant for $n < \omega$
\itemitem{ (b) }  if $u \triangleleft v^\ell \in [B]^{< \aleph_0}$ for
$\ell = 1,2$ then we can find $v_i \in [\lambda]^n$ for $i < \lambda,
u \subseteq v_i,\text{ min}(v_i \backslash u) \ge i$, and
$i < j \Rightarrow F(v^1 \cup v^2) = F(v_i \cup v_j)$.
\endroster}
\endroster
\endproclaim
\newpage

\nocite{Sh:276}
\nocite{Sh:288}
\nocite{Sh:289}
\nocite{Sh:473}
\nocite{Sh:546}

\bigskip
REFERENCES
\bigskip

\bibliographystyle{literal-plain}
\bibliography{lista,listb,listx,listf}

\newpage
\head {Private Appendix: old \S2} \endhead
\bigskip

We return here to consistency of statements of the form  $\chi \rightarrow  
[\theta]^2_{\sigma,2}$ (i.e. for every  $c:[\chi]^2 \rightarrow \sigma$  
there is  $A \in [\chi]^\theta$ such that on $[A]^{\le 2}$, $c$  has at 
most two values), (when  $2^{<\mu} \ge \chi > \theta^{>\mu}$,  of 
course).  In [Sh276,\S2] this was done for  $\mu = \aleph_0$, $\chi = 
2^\mu$,  $\theta = \aleph_1$ and  $\chi$  quite large (in the original 
universe it is an Erdos cardinal).  Originally, it was written for any
$\mu = \mu^{<\mu}$ ($\chi$  measurable in the original universe) but
because of the referee urging it is written up there for  $\mu = \aleph_0$
only; though with
an eye on the more general result which is only stated.  In [Sh288] the main 
objective is to replace colouring of pairs by colouring of $n$-tuples
(and even $(< \omega)$-bytes).  By [Sh284] we can lower the consistency
strength.  The main point of this section is to increase  $\theta$,
and this time do it for 
$\mu = \mu^{>\mu} > \mu$,  too.
\medskip

\noindent
\definition{2.1 Definition}  A forcing notion $Q$ satisfies
$\ast^\epsilon_\mu$ where $\epsilon$  is a limit ordinal $< \mu$, if
player I has a winning strategy in the following game  $\ast^\epsilon_\mu[Q]$
defined as follows:

\noindent
Playing: the play finishes after  $\epsilon$ moves.   

In the $\alpha$-th the move:   

Player I --- if  $\alpha \ne 0$  he chooses  
$\langle q^\alpha_\zeta:\zeta < \mu^+ \rangle$  such that  
$q^\alpha_\zeta \in Q$        

and  $(\forall \beta < \alpha)(\forall \zeta < \mu^+)p^\beta_\zeta \le  
q^\alpha_\zeta$ and he chooses a regressive function $f_\alpha:\mu^+ 
\rightarrow \mu^+$ (i.e.  $f_\alpha(i) < 1 + i)$;  if $\alpha = 0$  let  
$q^\alpha_\zeta = \emptyset_Q$, $f_\alpha = \emptyset$. \newline   

Player II --- he chooses  $\langle p^\alpha_\zeta:\zeta < \mu^+ \rangle$  
such that  $q^\alpha_\zeta \le p^\alpha_\zeta \in Q$.

\noindent
The Outcome:  Player I wins provided whenever $\mu < \zeta < \xi < \mu^+$,
$cf(\zeta) = cf(\xi) = \mu$  and  $\dsize \bigwedge_{\beta < \epsilon}
f_\beta(\zeta) = f_\beta(\xi)$  the set  
$\{p^\alpha_\zeta:\alpha < \epsilon \} \cup \{p^\alpha_\xi:\alpha < 
\epsilon \}$  has an upper bound in  $Q$.
\enddefinition
\bigskip

\remark{Remark}  In the case $\mu = \aleph_0$ we can use Knaster condition
instead.
\endremark
\bigskip

\definition{2.2 Definition}  Assume  $P$, $R$ are forcing notions,
$P \subseteq  R$,  $P \lessdot R$.  Moreover, for simplicity, for $r \in R$
a member  $r \restriction P \in  P$  is defined such that:
\medskip
\roster
\item "{(a)}"  $r \restriction P \le r$,
\item "{(b)}"  $(\forall p)(r \restriction P \le p \in P \Rightarrow r)$,
$p$  compatible in  $R$
\item "{(c)}"  $r^1 \le r^2 \Rightarrow  r^1 \restriction P \le r^2
\restriction P$
\endroster
\medskip 

\noindent
(so  $P \lessdot R$). 

Let  $(R,P)$  satisfies  $\ast^\epsilon_\mu$ mean:
\medskip
\roster
\item "{$(\alpha)$}"  $P$  satisfies  $\ast^\epsilon_\mu$
\item "{$(\beta)$}"  If  $St_1$ is a winning strategy for player
I in the game  $\ast^\epsilon_\mu[P]$  then in the following game called  
$\ast^\epsilon_\mu[P,R;St_1]$ the first player has a winning strategy:
\endroster
\medskip    

\noindent
Playing: as before but  $\langle <q^\theta_\zeta \restriction P:\zeta < 
\mu^+>,<p^\theta_\zeta \restriction P:\zeta < \mu^+>:\alpha < \epsilon
\rangle$
is required to be a play of $\ast^\epsilon_\mu[P]$  in which first player 
uses the strategy  $St_1$.
\medskip

\noindent
The outcome:  Player I wins provided that: if  $\zeta < \xi < \mu^+$,  
$cf(\zeta) = cf(\xi) = \mu$,  $\dsize \bigwedge_{\beta < \alpha} f_\beta
(\zeta) = f_\beta(\xi)$  and  $r^\ast \in P$  is a $\le_P$-upper bound of
$\{p^\alpha_\zeta \restriction P:\alpha < \epsilon \} \cup \{p^\alpha_\xi
\restriction P:\alpha < \epsilon \}$,  then  $\{r \} \cup \{p^\alpha_\zeta:
\alpha < \epsilon \} \cup  \{p^\alpha_\xi:\alpha < \epsilon \}$  has an
upper bound in  $R$.
\enddefinition
\bigskip

\definition{2.3 Definition/Lemma}  Let  $\mu = \mu^{<\kappa} < \kappa \sigma
 = cf\kappa \le \lambda \le \chi$.  (Usually fixed hence suppressed in the 
notion).  We define and prove the following by induction on (the ordinal)  
$\alpha$: \newline 
1) [Def]  We let  ${\Cal K}^{\epsilon,\alpha} = 
{\Cal K}^{\epsilon,\alpha}_{\mu,\kappa,\lambda,\chi}$ is the family of 
sequences  $\bar Q = \langle P_\beta,{\underset\tilde {}\to Q_\beta},
a_\beta,\bold t_\beta:\beta < \alpha \rangle$  such that:
\medskip
\roster
\item "{(a)}"  $\langle P_\beta,{\underset\tilde {}\to Q_\beta}:\beta <
 \alpha \rangle$  is a ($< \mu$)-support iteration (so $P_\alpha =
\text{ lim}\bar Q$  denotes the natural limit)
\item "{(b)}"  $a_\beta \subseteq \beta$,  $|a_\beta| < \kappa$, $[\gamma
\in a_\beta \Rightarrow  a_\gamma \subseteq a_\beta]$
\item "{(c)}"  ${\underset\tilde {}\to Q_\beta}$ has cardinality $< \lambda$
and is a $P_{a_\beta}$-name (see part $(*)$ below)
\item "{(d)}"  $\bold t_\beta$ is a truth value; i.e. true or false
\item "{(e)}"  if  $\bold t_\beta$ is true, $\beta < \gamma \le \alpha$
then  $(P_\gamma,P_\beta)$  satisfies  $\ast^\epsilon_\mu$.
\endroster
\medskip 

\noindent
2) [Def]  For  $\bar Q$  as above:
\medskip
\roster
\item "{(a)}"  $a \subseteq \alpha$  is called $\bar Q$-closed if
$\beta  \in  a \Rightarrow  a_\beta \subseteq  a$
\item "{(b)}"  for a $\bar Q$-closed subset $a$  of $\alpha$  we let    
\newline
$P_a = \biggl\{  p \in  P_\alpha:\text{Dom}(p) \subseteq  a$ and
for each $\beta  \in  \text{ Dom}(p)$ we have $p(\beta)$ \newline
$\qquad \qquad$ is a $P_{a \cap \beta}$-name $\biggr\}$ \newline  
$P^\ast_a = \left\lbrace p \in  P_\alpha:\text{Dom} p \subseteq  a\right.$ and 
for each  $\beta \in \text{ Dom}(p)$ we have  $p(\beta)$ is a
$P_{a_\beta}$-name \newline
$\qquad \qquad$ and if ${\underset\tilde {}\to Q_\beta} \subseteq V$
then  $p(\beta)$  is from  $V$  (see 0.B.) $\biggr\}$.
\endroster
\medskip
\noindent
3) [Lemma]  For  $\bar Q$  as above,  $\beta < \alpha$
\medskip
\roster
\item "{(a)}"  $\bar Q \restriction \beta \in  K^{\epsilon,\beta}$
\item "{(b)}"  if  $a \subseteq  \beta$  then $a$ is $Q$-closed iff
$a$  is ($\bar Q \restriction \beta$)-closed
\item "{(c)}"  if  $a \subseteq \alpha$  is $\bar Q$-closed then so is
$a \cap  \beta$.
\endroster
\medskip
\noindent
4) [Lemma]  For  $\bar Q$  as above, and  $p < \alpha$,  $P_\beta  
\lessdot P_\alpha$,  moreover, if  $p \in  P_\alpha$,  $p \restriction \beta  
\le q \in  P_\beta$ then  ($p \restriction [\beta \cup \alpha]) \cup  q \in  
P_\alpha$ is a lub of $p,q$. \newline
\noindent
5)< [Lemma]  For  $\bar Q$  as above
\medskip
\roster
\item "{(a)}" $P^\ast_\alpha$ is a dense subset of  $P_\alpha$ 
\item "{(b)}"  if  $a$  is $\bar Q$-closed then $P_a \lessdot P_\alpha$,  
$P^\ast_a$ is a dense subset of  $P_a$.
\endroster
\medskip
\noindent
6)  If  $a \subseteq  b$  is a $\bar Q$-closed set then  
$(P^\ast_b,P^\ast_a)$  satisfies  $\ast^\epsilon_\mu$. \newline
\noindent
7)  The sequence  $\bar Q = \langle P_\beta,{\underset\tilde {}\to Q_\beta},
a_\beta,\bold t_\beta:\beta < \alpha \rangle$  belongs to
${\Cal K}^{\epsilon,\mu}$ if $\alpha$  is a limit ordinal,
$\dsize \bigwedge_{\gamma < \alpha} \bar Q \restriction \gamma \in
\kappa^{\epsilon,\gamma}$ for some club  $C$  of  $\alpha$,
$\dsize \bigwedge_{\gamma \in c} \bold t_\gamma  =$ true. \newline
\noindent
8)  The sequence  $\bar Q = \langle P_\beta,{\underset\tilde {}\to Q_\beta},
a_\beta,\bold t_\beta:\beta < \alpha \rangle$  belongs to
${\Cal K}^{\epsilon,\mu}$ if $\alpha = \gamma + 1$, $a \subseteq \gamma$
is a ($\bar Q \restriction \gamma$)-closed set of cardinality  $< \kappa$,
$Q_\beta$ is a $P^\ast_{a_\beta}$-name of a forcing notion satisfying  
$\ast^\epsilon_\mu$ and $\bold t_\beta$  = true.
\enddefinition
\bigskip

\demo{Proof}  [Saharon].
\enddemo
\bigskip

\noindent
\proclaim{2.4 Theorem}  Suppose  $\mu = \mu^{<\mu} < \kappa = \lambda < \chi$,
$\chi$  is measurable.  For some forcing notion  $P$  of cardinality  $\chi $,
$\mu$-complete not collapsing cardinalities not changing cofinalities we have:
\newline 
\medskip
$\Vdash_P ``2^\mu = \chi$  and for every  $\sigma < \chi$  and $\theta < 
\kappa$  we have $\chi \rightarrow  [\theta]^2_{\sigma,2}$"  (and Axiom as
in [Sh280]) (and the parallel of  $M$  for forcing notions satisfying  
$\ast^\epsilon_\mu$).
\endproclaim
\bigskip

\demo{Proof}  Fix  $\epsilon < \mu$,  ${\Cal K}^\alpha  = 
{\Cal K}^{\epsilon,\alpha}$,  ${\Cal K} = \dsize \bigcup_{\alpha < \chi}
{\Cal K}^\alpha$.  By preliminary forcing without loss of generality  
``$\chi$ measurable" is preserved by forcing with  $({}^{\chi >} 2,
\trianglelefteq)$,($=$ adding a Cohen subset of  $\chi$).  Let us
define a forcing notion  $R$:
\medskip
\noindent
$R = \{\bar Q:\bar Q \in {\Cal K}$  for some  $\alpha < \chi$  and  $\bar Q 
\in  H(\chi)\}$  \newline
ordered by:  $\bar Q^1 \subseteq  \bar Q^2$ iff  $\bar Q^1 = 
\bar Q^2 \restriction \ell g(Q^1)$. \newline
\noindent
As  $R$  is equivalent to  $({}^{\chi >} 2,\triangleleft)$ we know that in
$V^R$,  $\chi$  is still measurable.  Let  $\bar Q^\oplus  = 
\langle P_\beta,{\underset\tilde {}\to Q_\beta},a_\beta,
\bold t_\beta:\beta < \chi \rangle$  be $\cup G_R$ and  $P_\chi$ be the
limit so  $P^\ast_\chi \subseteq  P_\chi$ is a dense subset.  Now
$R * {\underset\tilde {}\to P^\ast}$ is the forcing  $P$  we have promised.
The non-obvious point is  
$\Vdash_{R * P^\ast} ``\chi * [\theta]^2_{\sigma,2}$ (where  
$\theta < \chi,\sigma < \chi$).  So suppose  $(r^\ast,p) \in  R * 
{\underset\tilde {}\to P^\ast}$,  $(r^\ast,p) \Vdash$ ``the colouring
${\underset\tilde {}\to \tau}:[\chi]^2 \rightarrow  \sigma$  is a
counterexample".  Let  $\chi_1 = \left( 2^\chi \right)^+$.  Let
$G_R \subseteq  R$  be generic over  $V$,  
$r^\ast \in G_R$.  By [Sh289],xx in  $V^R$ we can find  $B \in  
[\chi]^\chi$ and  $\langle M_s:s \in [B]^{< \aleph_0} \rangle$  and  
$\langle f_{s,t}:(s,t) \in  \dsize \bigcup_{n < \omega}
[{\underset\tilde {}\to B}]^n \times [{\underset\tilde {}\to B}]^n \rangle$
such that  $M_s \prec (H(\chi_1)^{V[G_R]},H(\chi_1),\in)$,  
$\{\chi,G_R,p,{\underset\tilde {}\to \tau} \} \in  M_s$,
$M_s \cap  M_t = M_{s \cap t}$,  
$[M_s]^{\le(\kappa + \lambda + 2^\mu)} \subseteq M_s$,
$f_{s,t}$ isomorphism from  $M_s$ onto  $M_t$ as there. 

Let  $C = \biggl\{ \delta < \chi:\delta  = \text{ sup}(\beta \cap  
\delta)$ and  $(s \in  [\delta]^{< \aleph_0} \Rightarrow  M_s \cap  
\chi \subseteq  \delta) \biggr\}$. \newline
\medskip
\noindent
Let  $\gamma(\ast) = \text{ Min}(B)$.  Now for  $p \in P^\oplus \cap
M_{\gamma(\ast)}$ and  $u = \{c_1,c_2\} \in [\sigma]^{\le 2}$ let us
define the statement
\medskip
$(\ast)^u_p:$ if  $p \le p^0 \in  p^\ast \cap  M_{\gamma(\ast)}$ then we
can find  $p^1,p^2 \in  P^\ast \cap  M_{\gamma(\ast)}$, $p^0 \le  
p^1_\ast$,  $p^0 \le p^2$ such that:    

for  $\gamma_1 \ne \gamma_2$, $\gamma_2 \in B$, $\gamma_2 \in B$, we can
find  $r_1$, $r_2 \in  P^\ast  \cap M_{\{\gamma_1,\gamma_2\}}$
(so $\text{ Dom } r_\ell \subseteq  M_{\{\gamma_1,\gamma_2\}} \cap \chi$)
such that:

$$
r_\ell \Vdash ``{\underset\tilde {}\to \tau}(\{\gamma_1,\gamma_2\}) = c_\ell"
$$
 
$$
r_\ell \restriction \left( \lambda  \cap  M_{\{\gamma_\ell \}}\right)  = 
f_{\{\gamma(\ast)\},\{\gamma_\ell \}}(P^1)
$$
 
$$
r_\ell \restriction \left( \lambda \cap  M_{\{\gamma_{3-\ell}\}}\right)  = 
f_{\{\gamma(\ast)\},\{\gamma_{3-\ell}\}}(P^2).
$$
\medskip

\noindent
Easily  $I = \{ p \in  P^\ast  \cap  
M_{\{\gamma (\ast )\}}:$ \text{for some } $u$, $(\ast)^u_p$ 
hold $\}$ is a dense subset of  $P^\ast \cap  M_{\{\gamma(\ast)\}}$,
but this partial forcing satisfies the $\mu^+$-c.c. hence we can find
$\bold I^\ast = \{p_\zeta:\zeta < \mu \} \subseteq  I$  a maximal
antichain of  $P^\ast \cap M_{\{\gamma(\ast)\}}$ hence of $P^\ast$, and let 
\medskip
$(*)^{(c_1[p],c_2[p])}_\beta$ $\qquad$ for  $p \in \bold I^*$. 
\medskip

As  $G_R$ was any generic subset of  $R$  to which  $r^\ast $ belongs we have 
$R$-names ${\underset\tilde {}\to \gamma}(*) \langle
(p_\zeta,c_1(p_\zeta),c_2(p_\zeta)):\zeta < \mu \rangle$,
$\langle {\underset\tilde {}\to M_s}:s \in [{\underset\tilde {}\to B}]
^{<\aleph _0} \rangle$, $\langle {\underset\tilde {}\to f_{s,t}}:
(s,t) \in \dsize \bigcup_{n < \omega}[{\underset\tilde {}\to B}]^n \times
[{\underset\tilde {}\to B}]^n \rangle$  forced by  $r^*$ to be as above.
Without loss of generality 
$r^*$ force  a  values  $\gamma(*),M_\emptyset,M_{\{\gamma(*)\}}$,  
$\langle (p^\ast_\zeta,c^\ast_1(p_\zeta),c^\ast_2(p_\zeta)):\zeta  < 
\mu \rangle$. 

We now choose by induction on  $\zeta \le \sigma + 1$,  
$Q^\zeta,\alpha^\zeta,\gamma^\zeta$ such that:
\medskip
\roster
%\widestnumber\item{(A)(a)}
\item "{(A)(a)}"  $Q^\zeta \in  R$
\item "{(b)}"  $\bar Q^0 = r^\ast$
\item "{(c)}"  $\ell g(\bar Q^\zeta)  = \alpha^\zeta$
\item "{(d)}"  $\xi  < \zeta  \Rightarrow  \bar Q^\xi  = 
\bar Q^\zeta \restriction \alpha^\zeta$
\item {(e)}"  $\langle \alpha^\zeta:\zeta \le \sigma + 1 \rangle$
is strictly increasing continuous
\item "{(f)}"  $\alpha^\zeta  < \gamma^\zeta  < \alpha ^{\zeta +1}$ 
\item "{(g)}"  $\bar Q^{\zeta + 1} \Vdash_R ``\gamma^\zeta  \in  
{\underset\tilde {}\to B}$"
\item "{(h)}"  $\bar Q^{\zeta + 1}$ forces a value to  $\langle 
{\underset\tilde {}\to M_s}:s \in [{\underset\tilde {}\to B} \cap 
(\gamma^\zeta + 1)^{<\omega} \rangle  = \langle 
{\underset\tilde {}\to M_s}:s \in [{\underset\tilde {}\to B_\zeta}]
^{<\omega} \rangle$
\item "{(B)(a)}" if  $\zeta \le \sigma + 1,cf(\zeta) > \mu$ then:
{\roster
\itemitem{ (a) }  $a^{\bar Q^{\zeta +1}}_\zeta  = \cup \{ \chi  \cap  
M_{\{\xi _1,\xi _2\}}:\{\xi _1,\xi _2\} \in  [\{\gamma_\epsilon:\epsilon  < 
\zeta \}]^{\le 2} \}$
\itemitem{ (b) }  ${\underset\tilde {}\to Q^{\bar Q^{\zeta +1}}
_{\alpha_\zeta}} = \{ h:h$ a function   
$\text{ Dom}(h) \in [\{\gamma_\epsilon:\epsilon < \zeta \}]^{<\mu},
\text{ Rang }h \subseteq \mu\}$.        
\endroster}
\endroster
\medskip

\noindent
For  $\gamma_\epsilon \in \text{ Dom } h$,  
$f_{\{\gamma(*)\},\{\gamma_\epsilon \}}(p^*_{h(\gamma_\epsilon)}) 
\in {\underset\tilde {}\to G_{P_\zeta}}$ for $\epsilon_1 < \epsilon_2 
< \zeta$,  we can find  $p^1$, 
$p^2$, $r_1$, $r_2$ such that:  $p^\ast_\zeta \le p^\ell  \in  
M_{\{\eta(*)\}} \cap P^\ast_\chi$,        

$r_\ell  \in  P^\ast _\chi  \cap  
M_{\{\gamma _{\epsilon _1},\gamma _{\epsilon _2}\}},$        
\smallskip
$r_\ell \restriction (\lambda  \cap  M_{\{\gamma _{\epsilon _\ell }\}}) = 
f_{\{\gamma (\ast ),\{\gamma _{\epsilon _1}\}}(p^1),$        
\smallskip
$r_\ell \restriction (\lambda  \cap  M_{\{\gamma _{\epsilon _{3-\ell }}\}}) = 
f_{\{\gamma (\ast )\},\{\gamma _{\epsilon _{3-\ell }}\}}(p^2),$        
\smallskip
$r_1 \in {\underset\tilde {}\to G_{P_\zeta}}$ or  $r_2 \in  
{\underset\tilde {}\to G_{P_\zeta}}$. 
\medskip

The point is to verify that defining
$Q^{\bar Q^{\zeta +1}}_{\alpha _\zeta }$ 
as we do in clause $(B)(a)$ is allowable; i.e. that the condition concerning 
$*^\epsilon_\mu$ from Definition 2.2 holds.

\noindent
Saharon more.
\enddemo
\bigskip

\remark{2.5 Remark}
We can in 2.4 replace ``measurable", by $k^2_1$-Mahlo, but it is 
not straightforward; $.e.g$. we may use the following version of squared 
demand. 

Let  $\chi $  be $k$-Mahlo,  $\kappa  < \chi $ 
\medskip
$(\ast )^{\kappa ,k}_\kappa $ there is  $\bar A = \langle A_\alpha :\alpha  < 
\chi \rangle $  and  $\bar C = \langle C_\alpha :\alpha  < S\rangle $  such 
that:
\medskip
\roster
\item "{(a)}"  $S \subseteq  \{\delta < \chi:cf\delta \le \kappa \}$,
$\{\delta \in S:cf\delta = \kappa \}$  is a stationary subset of  $\chi$
\item "{(b)}"  $C_\alpha  \subseteq  \alpha$,  $[\beta \in c_\alpha
 \Rightarrow  c_\beta  = c_\alpha  \cap  \beta ]$,
otp\ $c_\alpha  \leq  \kappa $,  $c_\alpha $ a club of  $\alpha $    
\item "{(c)}"  $A_\alpha  \subseteq  \alpha ,$    
\item "{(d)}"  $\beta \in  c_\alpha  \Rightarrow  A_B = A_\alpha \cap \beta$
\item "{(e)}"  $\{ \lambda  < \chi :\lambda$  inaccessible, and for 
every  $x \subseteq \lambda  \{\alpha < \chi:cf\alpha = \kappa < 
\alpha ,x \cap  \alpha  = A_\alpha \}$  is stationary $\}$  
not only is stationary but is $k$-Mahlo.
\endroster
\medskip

\noindent
This can be obtained e.g. by iteration with Easton support, in which for each
inaccessible we add  $\bar A$, $\bar c$  satisfying (a)-(d) above, each 
condition being an initial segment. 
\bigskip

\centerline {$* \qquad  * \qquad  *$} 
\bigskip

We do now what we can for colouring of triples and more.

\noindent
By \cite{Sh:288},\S4 we can deal with $\mu$ supercompact and, of course, with  
$\mu = \aleph_0$.
\endremark
\newpage

\head {\S3} \endhead
\bigskip

\subhead {3.1 Content} \endsubhead   $\tau$  is a vocabulary $=$ set of
predicates including  $<, K^\tau$ the class of $\tau$-structures,
$\mu$  such that  $<^\mu$ is a 
(linear) well ordering,  $M$, $N \in  K^\tau$ are compatible if
$M \restriction (|M| \cap  |N|) = N \restriction (|M| \cap |N|)$, and
the  $M \cup  N$  is  $(|M| \cup  |N|,\dotsc,R^M \cup  R^N,\ldots)
_{R \in \tau }$.  Universe of a model  $M$  is  $|M|$  
so for  $M \in  \tau$,  $M = (|M|,\dotsc,R^M,\ldots)_{R\in \tau}$,
$M \in  K$  is standard is  $|M|$  is a set of order,  $<^M$ usual
ordering,  $\tau  = \tau \backslash \{<\},M^- = M \restriction \tau^-$.
\bigskip

\definition{3.2 Definition}  1) We call
$M \in  K^\tau  2$-connected if for no  $N_1N_2$ 
do we have  $N_1 \subseteq  M$,  $N_2 \subseteq M$, $M^- = N^-_1 \cup  N^-_2$
and  $\|N_1 \cap  N_2\| \le 1$.
\newline
2)  cont$^\kappa(M) = \{A \subseteq  |M|:A \subseteq  M$, $\text{otp}(A) <
\kappa $, $M \restriction A$  is 2-connected $\}$, \newline
cont$^\kappa  \cong  (M) = \{(M \restriction A) / \cong:A \subseteq
|M|,|A| < \kappa,M \restriction A$  is 2-connected $\}$.
\newline
3)  For  $D$  a set of members of  $K^\tau$,  $D / \cong$  is
$\{M/\cong :M \in D\}$. \newline
4)  If  $D$  is a set of 2-connected members of  $K^\tau$, $\sigma < \kappa 
< \infty $,  $K^\tau _\kappa [D] = \{M \in  K^\tau :$cont$_\kappa  \cong  (M) 
\subseteq  D/\cong \}$  \newline
(so without loss of generality  $(\ast )[M \in  N 
\Rightarrow  \|M\| < \kappa ][A \subseteq  M \in  D \and 
M \restriction A 2$-connected  $\Rightarrow  M \restriction A \in  D]$.
\newline
5)  $M \in  K_\tau$ is $m$-indecomposable if for every  $x$, $y \in  M$  there
is  $A$,  $\{x,y\} \subseteq  A \subseteq  M$,  $|A| < m$,  
$M \restriction A$ 2-connected. \newline
6) $id \text{ con}_{m,\kappa}(M) = \{A \subseteq  M:|A| < \kappa $ and
$M \restriction A$ is $m$-indecomposable $\}$
\newline
$\text{id com}_{m,\kappa} \cong |M| = \{(M \restriction A) / \cong:A \in  
\text{ id com}_{m,\kappa }(M)\}$. \newline
If  $\kappa > \|M\|$  we omit it.
\enddefinition
\bigskip

\proclaim{3.3 Theorem}  Suppose
\medskip
\roster
\item "{(a)}"  $\mu < \lambda,\mu = \mu^{<\mu},\lambda$ measurable,  
$\theta$
\item "{(b)}"  $1 < m < \omega,D$  a set of 2-connected members of  $K^\tau $ 
satisfying $(*)_m$ of 2-2, $|\tau|$
\item "{(c)}"  $M \in K^\tau_k[D],\|M\| < \mu$ and $|M|$ is $\alpha^* \le 
\mu$.
\endroster
\medskip 

\noindent
Then for some  $P$  and ${\underset\tilde {}\to N}$:
\medskip
\roster
\item "{(A)}"  $P$  is a $\mu$-complete $\lambda$-c.c. forcing notion 
of cardinality $\lambda$
\item "{(B)}"  ${\underset\tilde {}\to N}$ is a $P$-name of a 
standard member of $K^\tau_k[D]$  with universe  $\lambda$ 
\item "{(C)}"  ${\underset\tilde {}\to N}$ ``if $c:^{\omega >}\lambda \ast 
\mu$  then for some  $A \subseteq \lambda$: \newline
$\text{otp } A = \alpha^*,{\underset\tilde {}\to N} \restriction A \cong \mu$
and $(*)$ if  $n < \omega$, \newline
$\alpha_0 < \alpha_1 < \cdots < \alpha_{n-1} < \alpha_n$ are from $A$,
\newline       
${\underset\tilde {}\to N} \restriction \{\alpha_0,\alpha_1,\dotsc,
\alpha_{n-2},\alpha_{n-1}\} \cong {\underset\tilde {}\to N} \restriction 
\{\alpha_0,\alpha_1,\dotsc,\alpha_{n-2},\alpha_n\}$ then \newline
$c(\{\alpha_0,\dotsc,\alpha_{n-2},\alpha_{n-1}\}) = 
c(\{\alpha_0,\dotsc,\alpha_{n-2},\alpha_{n-1}\})$
\item "{(D)}"  $\Vdash_P ``\text{id cont}_m \cong ({\underset\tilde {}\to N})
\subseteq \text{ id cont}_m \cong (G)$"
\item "{(E)}"  $\Vdash_P ``\lambda = \mu^+"$.
\endroster
\endproclaim
\bigskip

\demo{Proof}:
\enddemo
\bigskip

\noindent
A Stage:  The forcing:
\medskip
\roster
\item "{$(\alpha)$}"  Let  $Q$  be the set of $p = (s^p,N^p,\varphi^p) = 
(s,N,\varphi)$  such that
{\roster
\itemitem{ (a) }  $s \subseteq \lambda$, $(s) < \mu$
\itemitem{ (b) }  $N^p \in K^\tau_m[D]$  has universe ${}^k s = 
\{\langle \alpha_0,\dotsc,\alpha_{m-1} \rangle:\alpha_0 > \alpha_1 > 
\cdots \alpha_{m-1}\}$
\itemitem{ (c) }  $<^{N^p}$ is the lexicographic order on ${}^k s$
\itemitem{ (d) }  for  $R \in \tau^-$ if  
$\langle \bar \alpha^0,\dotsc,\bar \alpha^{m-1} \rangle \in R^{N^P}$ then 
for  
$\ell \ne k(< n) \bar \alpha^\ell \sim \bar \alpha^m$ which means  
$\alpha^\ell_0 > \alpha^m_0 > \alpha^\ell_0 > \alpha^m_1 > 
\cdots$ or $\alpha^m_0 > \alpha^\ell_0 > \alpha^m_1 > 
\alpha^\ell_1 > \cdots$.  We shall not strictly distinguish 
between  $N^p$ and  $N^p \restriction \tau^-$
\itemitem{ (e) }  $\varphi$  is a function,$\text{Dom } \varphi = 
\text{ id cont}_{m,\kappa}(M)$  and 
$\varphi(N \restriction A)$  is an embedding of  $N \restriction A$ 
into  $M$
\itemitem{ (f) }  if  $A,B \subseteq \text{ Dom } \varphi,A \subseteq B$
moreover $A$  is an initial segment of  $B$ by  $<^N$ then 
$\varphi (N \restriction A) \subseteq  \varphi (N \restriction B)$
\itemitem{ (g) }  for  $\alpha < \lambda,p \restriction \alpha \le p$  
where  $p \restriction \alpha = (s^p \cap \alpha,N^p \restriction {}^m(s^p 
\cap \alpha),\varphi^p \restriction \text{ id cont}_m(N^p \restriction
{}^m(s^p \cap \alpha))$, and $\le$ is defined in $(\beta)$ below.
\endroster}
\item "{$(\beta)$}"  The order on  $Q:p \le q$ iff:
{\roster
\itemitem{ (a) } $s^p \subseteq  s^q$         
\itemitem{ (b) } $N^p \subseteq  N^q$         
\itemitem{ (c) } $\varphi^p \subseteq  \varphi ^q$         
\itemitem{ (d) } if  $A \subseteq \text{ cont}_m(N^q),|A \cap N^p| > 1$
then  $A \subseteq N^p$         
\itemitem{ (e) } if  $A \subseteq N^q,|A| < m$, $|A \cap N^p| > 1$,
$A \subseteq N^p$ then we can find  $A_1$, $A_2$ such that $A = A_1 \cup A_2$,
$|A_1 \cap  A_2| \le 1$,  $A_1 \cap  A_2 \subseteq  N^p$,  
$N^p \restriction A \restriction \tau^- = N^p \restriction A_1 \restriction
\tau^- \cup N^p \restriction A_2 \restriction \tau^-$         
\itemitem{ (f) } if  $B \in \text{ Dom } \varphi^q$, $A = B \cap N^p \in
\text{ Dom } \varphi^p$ then  $B$  is an end extension of $A$.
\endroster}
\endroster
\bigskip  

Lastly,  $P = Q \times \text{ Levy}(\mu_2 < \lambda)$. But to simplify 
presentation we prove only that  $Q$  satisfies (A)-(D) (waiving (E)).
\medskip

\noindent
$B$ Stage:  Basic properties of the forcing notion:
\medskip
\roster
\item "{$(\alpha)$}"  if $p \le q$  in  $Q$, $A \subseteq N^q$, $|A| < m$,
then we can find  $A_\ell(\ell < \kappa)$  such that: \newline
``$\langle A_\ell \backslash N^p:\ell < \kappa \rangle$  is a partition of  
$A \backslash N^p$  \newline 
$|A_\ell \cap  N^p| \le  1$ \newline
$N^q \restriction A = N^q \restriction A_0 \restriction \tau^- + N^q 
\restriction A_1 \restriction \tau^- + 
\cdots + N^q \restriction A_{k-1} \restriction \tau^- + N^q \restriction 
(A \cap N^p)$
\newline
[why?  by induction on  $|A|$  using (A)(B)(e)].
\item "{$(\beta)$}"  if  $p \le q$  in  $Q$,  $B \in \text{ id cont}_m(q)$
then  $B \cap  N^p \in \text{ id com}_m(q)$ \newline    
[why?  if  $x,y \in B \cap N^p$ there is  $A \subseteq  B$,  $|A| < m$,  
$\{x,y\} \subseteq  A$,  $M \restriction A \, 2$-connected, by $(\alpha)$ 
above  necessarily  $A \subseteq  N^p$,  so we are done].
\item "{$(\gamma)$}"  $(\theta,\le)$  is transitive. \newline    
[why?  check $A(\beta)(a)$, $(b)$, $(c)$, $(d)$.  As for $(A)(\beta)(e)$ use 
$(\alpha)$ above and as for $(A)((\beta)(f)$ use $(\beta)$ above].
\item "{$(\delta)$}"  Assume  $q \in Q$,  $\alpha < \beta < \lambda$.  Then  
$(p \restriction \beta) \restriction \alpha = p \restriction \alpha \le 
p \restriction \beta \le p$. \newline    
[why?  by $(A)(\alpha)(g)$]. 
\item "{$(\epsilon)$}" if  $p \le q$  and  $\alpha < \lambda$  then  
$p \restriction \alpha \le q \restriction \alpha$. \newline    
[why?  check].
\item "{$(\zeta)$}"  Assume  $p_\ell \in Q(\ell = 0,1,2) p_0 \le p_1$, $p_0 
\le p_2$,  $s^{p_0} = s^{p_1} \cap  s^{p_2}$ and  $q = (s^q,N^q,\varphi^q)$,
$s^q = s^p \cup  s^{p_2}$,  $N^r \restriction \tau^- = N^{p_1} \restriction
\tau^- \cup  
N^{p_2} \restriction \tau^-$ and  $\varphi^{p_1} = \varphi^{p_1} \cup  
\varphi^{p_2}$. \newline
Then  $p_\ell \le q \in Q$.
\newline    
[why?  we prove it by induction on  $\text{sup}(s^q)]$.
\endroster
\medskip

\noindent
First check  $q \in Q$;  i.e.  $(A)(\alpha)$, now $(A)(\alpha)(a)$, $(b)$,
$(c)$, $(d)$ should be clear.  For $(A)(\alpha)(d)$ we need: 
\medskip
\roster
\item "{$(*)_0$}" $\text{id cont}_m(N^q) = \text{ id cont}_m(N^{p_1}) 
\cup \text{ id cont}_m(N^{p_2})$.
\endroster
\medskip

\noindent
For this it suffices to: 
\roster
\medskip
\item "{$(*)_2$}"  if  $A \in \text{ id cont}_m(N^r)$ then 
$A \subseteq  N^{p_1}$ or  $A \subseteq  N^{p_1}$ which holds as if 
$x \in  A \backslash N^{p_1}$,  $y \in  A \backslash N^{p_2}$, there 
is  $B$,  $\{x,y\} \in  B \subseteq  A$,  $|B| < 
m$,  $N^r \restriction B$  is 2-connected, but by $(\alpha)$ above 
applied twice we get contradiction.
\endroster
\medskip

\noindent
We have shown 
\medskip
\roster
\item "{$(*)_3$}"  if  $A \in \text{ cont}_m(N^r)$, then $A \subseteq 
N^{p_1}$ or  $A \subseteq  N^{p_2}$.
\endroster
\medskip

\noindent
Now $(A)(\alpha)(f)$ is easy, as for $(A)(\alpha)(g):(A)(\beta)(a)$, $(b)$,
$(c)$ are easy.  $(A)(B)(d)$ holds by $(*)_3$ above.  As for 
$(A)(\beta )(e)$, let  $B \subseteq  N^q$,  $|B| < m$.  By $(\alpha)$ above 
applied to  $p_1 \restriction \alpha \le p_1$ we can find  
$\langle x_\ell,A_\ell:\ell < \ell^\ast \rangle$,  $N^{p_1} \restriction
A_\ell$ 2-connected,  $A_\ell  \cap  N^{p_1 \restriction \alpha } = 
\{x_\ell \}$,  
$N^{p_1}\kappa A = \sum_\ell N^{p_1} \restriction A_\ell + N^{p_1}
\restriction \left( A \cap  N^{p_1 \restriction \alpha }\right)$.
\bigskip 

\noindent
For each  $\ell$  if  $A_\ell \cap N^{p_0}....$ 
\medskip
Those notes were intended to give 
consistency results of the form.  We have  $k < \omega$  a class $K$ of 1 
model with set of finite submodels pregiven 
(closed under edgeless amalgamation) and set of forbidden large $k$-connected 
in appropriate sense submodels, we can have  $\forall \kappa \forall n
\forall M \in  \kappa \exists N \in  K$,  $M \rightarrow  (N)^n_\kappa$. 
\medskip

But the conditions on the set of finite submodels do not match the forcing 
construction such that appropriate amalgamations hold.  For ``no odd circle
of small length" fine should be rephrased.
\newpage

\head {Private Appendix} \endhead
\bigskip

(17.4.94):
1)  Replace $\aleph_0$ by higher...?
2)  Replace Ramsey by weakly compact??
\medskip
\noindent

%SAHARON AND ANDRZEJ: 
% still need to clean up
% ALICE

(4.5.94:2)  Remember \cite{Sh:481},\S2 make sense: 
we use the parallel of squared
scales - have start lecture and and finish correction \S1, correcting \S2.
Enough for $2^\mu \rightarrow [\mu^{++}]^2_3$. \newline
Write down:
\medskip
\roster
\item "{(a)}"  Mahlo
(use Easton support instead $R$, (which use $(< \chi)$-support, for a larger
cardinal) (generic from system $\langle N_s:i \rangle$?? preserve or only
continuity (see (c)(b);
\item "{(b)}"  check colouring
\item "{(c)}"  Try:
{\roster
\itemitem{ $(\alpha)$ }  some $\chi' < \chi,\chi' \rightarrow [\mu^+]^2_3$,
and $\chi' < \mu^{+ \omega}$ (so here the guessing of
$\langle M_{\gamma_\varepsilon}:\varepsilon \le \theta \rangle$ is not by
initial segments but by $\langle M_{\gamma_\varepsilon} \cap \delta_i:
\varepsilon \le \theta \rangle$ or
\itemitem{ $(\beta)$ }  back to what I thought last week, trying to force
almost disjoint set [Baumgartner]
\endroster}
\item "{(d)}"  289,\S5 in ZFC?
\endroster
\medskip

\noindent
5.5.94  Will check how preliminary forcing affect the systems of
\cite{Sh:289}, I though (April) will be like tree, so enough, have 
to check) the
proof of 481,\S1 indicates the Knaster in [Sh288] can be weakened sometimes.
\medskip

\noindent
11.5.94  In 1.7; the April idea of forcing almost disjoint sets may help to
$2^{\aleph_0}$-c.c. $\times \aleph_2$-c.c. $\Rightarrow 2^{\aleph_0}$-c.c.
\newline
The systems we need sometimes are just names for
${\underset\tilde {}\to c}(s)$: this is less: when it is not sufficient?
Different numbers.
\newpage

\head {Assignments from '92} \endhead
\bigskip

\noindent
0) explain for  $\mu = \aleph_0$ in \S2

\noindent
1) 2-3 fil

\noindent
2) 2-4 fil

\noindent
3) 2-5 fil

\noindent
4) 2-6 fil

\noindent
5) explicit c.c.c. condition --- easier if we represent  $\beta (i)$  as a 
set of ordinals 

(5A) $2^{\aleph_0} < \aleph_\omega$ 

(5B) $\chi \rightarrow [\theta]^2_{\sigma,2}$ for every
$\theta$;  $\sigma < \chi $ 

\noindent
6) $\chi  < 2^\mu ??$

\noindent
7) The reference of  $\lambda  \rightarrow  [\mu ]^3_\mu $ if  $2^{<\mu } < 
\lambda  \leq  2^\mu $ or\ so

\noindent
8) Concerning {\S}1, if  $2^{\aleph _0}$ is successor of regular?

\noindent
9) Can you guess square domain of  $\langle s_u:u \in  
[\mu ^+]^{\leq 2}\rangle   |s_4| \subseteq  \mu $  with  $\chi  < 
\mu ^{+\omega }??$

\enddocument
\bye